%% file: main_elsarticle-template-num.tex
\documentclass[preprint,11pt]{elsarticle}
\usepackage[showframe = false, headheight = 2cm, margin = 2cm, bottom = 3cm]{geometry}

\usepackage{amssymb}
\usepackage{bigints}
\usepackage{mathrsfs}
\usepackage{calc}
\usepackage{tikz}
\usepackage{xstring}
\usepackage{tkz-euclide}
\usepackage{graphics}
\usepackage{multicol}
\usepackage{adjustbox}
\usepackage{subcaption}
\usepackage{multirow}
\usepackage{etoolbox}
\usepackage{datetime}
\usepackage{amsmath, accents}
\usepackage{amsthm}
\usepackage{amssymb}
\usepackage{afterpage}
\usepackage{acro}
\usepackage{float}
\usepackage{relsize}
\usepackage{array}
\usepackage{textcomp}
\usepackage{comment}
\usepackage{diagbox}
\usepackage{listings}
\usepackage{tabularx}
\usepackage{ragged2e}
\usepackage{lscape}
\usepackage{enumitem}
\usepackage{eurosym}
\usepackage{footmisc}
\usepackage{todonotes}
\usepackage{mwe}
\usepackage{gensymb}
\usepackage{wasysym}
\usepackage{graphicx}
\usepackage{blindtext}
\usepackage{appendix}
\usepackage{multicol}
\usepackage{hyperref}
\usepackage{cleveref}
\usepackage{cancel}
\usepackage{bbm}
\usepackage{mathtools}
\usepackage{nomencl}
\usepackage{accents}
\usepackage{calligra}
\usepackage{miama}
\usepackage{empheq}
\usepackage{xcolor}
\usepackage{soul}
\usepackage{ulem}
\definecolor{codegreen}{rgb}{0,0.6,0}

\journal{}

\theoremstyle{definition}
\newtheorem{definition}{Definition}[section]
\DeclareMathOperator*{\argmin}{arg\,min}

\theoremstyle{remark}
\newtheorem{remark}{Remark}[section]

\begin{document}

\input{macros}
\usetikzlibrary{math}
\begin{frontmatter}



\title{Construction and application of an algebraic dual basis and the Fine-Scale Greens' Function for computing projections and reconstructing unresolved scales}

\author[1,3]{Suyash Shrestha}
\ead{s.shrestha-1@tudelft.nl; s.shrestha@upm.es}
\author[2,3]{Joey Dekker}
\author[3]{Marc Gerritsma}
\author[3]{Steven Hulshoff}
\author[4]{Ido Akkerman}
\affiliation[1]{organization={ETSIAE-UPM-School of Aeronautics, Universidad Politécnica de Madrid},
            addressline={Plaza Cardenal Cisneros 3},
            city={Madrid},
            postcode={E-28040},
            country={Spain}}
\affiliation[2]{organization={Delft University of Technology, Delft Institute of Applied Mathematics},
            addressline={Mekelweg 4},
            city={Delft},
            postcode={2628 CD},
            country={The Netherlands}}
\affiliation[3]{organization={Delft University of Technology, Faculty of Aerospace Engineering},
            addressline={Kluyverweg 1},
            city={Delft},
            postcode={2629 HS},
            country={The Netherlands}}
\affiliation[4]{organization={Delft University of Technology, Faculty of Mechanical Engineering},
            addressline={Mekelweg 2},
            city={Delft},
            postcode={2628 CD},
            country={The Netherlands}}



\begin{abstract}
\input{sections/01-abstract}
\end{abstract}



\begin{keyword}
duality; projection; variational multiscale; (Fine-Scale) Greens' function; Poisson equation; advection-diffusion equation
\end{keyword}

\end{frontmatter}

\tableofcontents


\input{sections/02-introduction}
\input{sections/03-proir_work}
\input{sections/04-dual_basis}
\input{sections/05-greens_function}
\input{sections/06-summary}


\bibliographystyle{elsarticle-num} 
\bibliography{references}

\newpage
\appendix
\input{sections/08-appendix}




\end{document}

%% file: macros.tex
\newcommand{\parddx}[2]{\dfrac{\partial #1}{\partial #2}}
\newcommand{\parddxdx}[2]{\dfrac{\partial^2 #1}{\partial #2^2}}
\newcommand{\E}[2]{\mathbb{E}^{#1, #2}}
\newcommand{\tE}[2]{\widetilde{\mathbb{E}}^{#1, #2}}
\newcommand{\tETranspose}[2]{\widetilde{\mathbb{E}}^{#1, #2^{T}}}
\newcommand{\ETranspose}[2]{\mathbb{E}^{#1, #2^{T}}}
\newcommand{\MassM}[1]{\mathbb{M}^{(#1)}}
\newcommand{\MassMP}[2]{\mathbb{M}^{(#1, #2)}}
\newcommand{\MassP}[3]{\mathbb{M}^{(#1, #2)}_{#3}}
\newcommand{\MassMinv}[1]{\mathbb{M}^{(#1)^{-1}}}
\newcommand{\weakForm}[2]{\int_{\Omega} #1 \wedge \star \: #2}
\newcommand{\boldCal}[1]{\boldsymbol{\mathcal{#1}}}

%% file: sections/01-abstract.tex
In this paper, we build on the work of [T. Hughes, G. Sangalli, VARIATIONAL MULTISCALE ANALYSIS: THE FINE-SCALE GREENS' FUNCTION, PROJECTION, OPTIMIZATION, LOCALIZATION, AND STABILIZED METHODS, SIAM Journal of Numerical Analysis, 45(2), 2007] dealing with the explicit computation of the Fine-Scale Green's function. The original approach chooses a set of functionals associated with a projector to compute the Fine-Scale Green's function. The construction of these functionals, however, does not generalise to arbitrary projections, higher dimensions, or Spectral Element methods.

We propose to generalise the construction of the required functionals by using dual functions. These dual functions can be directly derived from the chosen projector and are explicitly computable. We show how to find the dual functions for both the $L^2$ and the $H^1_0$ projections. We then go on to demonstrate that the Fine-Scale Green's functions constructed with the dual basis functions consistently reproduce the unresolved scales removed by the projector. 

The methodology is tested using one-dimensional Poisson and advection-diffusion problems, as well as a two-dimensional Poisson problem. We present the computed components of the Fine-Scale Green's function, and the Fine-Scale Green's function itself. These results show that the method works for arbitrary projections, in arbitrary dimensions. Moreover, the methodology can be applied to any Finite/Spectral Element or Isogeometric framework.

%% file: sections/02-introduction.tex
\section{Introduction}\label{sec:intro}
The origin of the Variational Multiscale (VMS) method stems from the series of papers on stabilisation techniques for finite-element methods found in \cite{HUGHES1986223, HUGHES1986341, HUGHES1986305, HUGHES1986329, HUGHES198685, HUGHES198797, HUGHES198785, HUGHES1989173, SHAKIB199135, SHAKIB1991141} along with \cite{shakib_1988}. The VMS methodology itself was first introduced in \cite{hughes_1995} as a re-interpretation of the SUPG formulation where it was shown that stabilisation methods could be derived from a solid theoretical foundation. Fundamentally, the multiscale framework entails incorporating the missing unresolved/fine-scale effects into the numerically computed large/resolved scales such that the numerical solution becomes the chosen projection of the exact solution. Moreover, it has served as a fundamental basis for the development of stabilised methods for applications such as turbulence modelling in the context of Large Eddy Simulation (LES), see \cite{munts_hulshoff_de_borst_2004, bazilevs_calo_zhang_hughes_2006, bazilevs_calo_cottrell_hughes_reali_scovazzi_2007, Holmen_Hughes_Oberai_Wells_2004, Hughes_Mazzei_Oberai_Wray_2001}. Other examples include \cite{stoter_turteltaub_hulshoff_schillinger_2018, bazilevs2006isogeometric, koobus_farhat_2004, levasseur_sagaut_chalot_davroux_2006, ten_eikelder_akkerman_2018} where VMS is employed for performing LES in Discontinuous Galerkin and Isogemetric frameworks, along with \cite{ten_eikelder_bazilevs_akkerman_2020} wherein Dicontinuity Capturing is tackled using multiscale theory. 

A key ingredient for formulating the VMS approach is the so-called \textit{Fine-Scale Greens' function}, which was first introduced in \cite{Hughes_Feijóo_Mazzei_Quincy_1998} and was formally characterised in \cite{hughes_sangalli_2007}. More specifically, in \cite{hughes_sangalli_2007} the explicit construction of the Fine-Scale Greens' function is demonstrated using the classic Greens' function and the choice of a projector. Although a general derivation of the Fine-Scale Greens' function is presented the actual construction for the Fine-Scale Greens' function for arbitrary projections in arbitrary dimensions is not entirely generalised. This is because the construction of the functionals associated with the projector, referred to as $\boldsymbol{\mu}$'s in \cite{hughes_sangalli_2007}, is not clarified for general cases. 

In the present work, we seek to build on \cite{hughes_sangalli_2007} by proposing a set of $\boldsymbol{\mu}$'s which is unique for every projector and generalises to arbitrary dimensions. To achieve this, we propose using a new set of basis functions, namely dual basis functions, to act as the $\boldsymbol{\mu}$'s. This enables the explicit construction of the Fine-Scale Greens' function for any given projector in arbitrary dimensions using the classic Greens' function. To start off, we briefly present a review of the work from \cite{hughes_sangalli_2007} in \Cref{sec:hughes_sangalli}. Subsequently, we present a discussion on the Mimetic Spectral Element Method (MSEM) and the associated primal and dual function spaces in \Cref{sec:dual_basis}. In this section, we also highlight the link between the chosen projector and the dual basis functions. Note that while we work with the MSEM throughout this study, the concepts are by no means limited to the MSEM. The presented work can be extended to Finite/Spectral Element or Isogeometric approaches by simply evoking the concept of duality introduced in this paper. 

In \Cref{sec:greens} we move on to construct the Fine-Scale Greens' function associated with the $H_0^1$ and $L^2$ projections for the 1D Poisson equation. Here we perform several numerical tests to demonstrate that the constructed Fine-Scale Greens' functions can exactly reconstruct all the fine scales truncated by the projection. Furthermore, we also show how the Fine-Scale Greens' function for the Poisson problem may be used in an advection-diffusion problem to reconstruct the corresponding fine scales. Although this involves solving coupled sets of equations instead of a direct computation, it overcomes the complexity of dealing with the classical Greens' function of the advection-diffusion problem. We close in \Cref{sec:proj_2d} by presenting the methodology applied to the 2D Poisson problem to demonstrate the generalisation to higher dimensions. 

%% file: sections/03-proir_work.tex
\section{Hughes and Sangalli approach}
\label{sec:hughes_sangalli}
The theory of variational multiscale analysis that this work builds upon was presented by Hughes and Sangalli in \cite{hughes_sangalli_2007}. This section serves as a brief review of the original work of \cite{hughes_sangalli_2007} where we lay down the fundamental concepts for formulating the VMS approach. First, the continuous setup is reviewed, followed by the finite-dimensional setting.

\subsection{Continuous Case}
In \cite[\S 2]{hughes_sangalli_2007}, the following problem is considered. Let $V$ be a Hilbert space with a norm $||\cdot||_V$ and a scalar product $(\cdot,\cdot)_V$. Let $V^*$ be the dual of $V$, with $\prescript{}{V^*}{\langle\cdot,\cdot\rangle}_{V}$ the duality pairing between $V^*$ and $V$. Let $\mathcal{L}: V \to V^*$ be a linear isomorphism. The problem can then be stated as: given $f \in V^*$, find $u \in V$ such that
\begin{equation}
    \mathcal{L} u = f.
    \label{eq:pde}
\end{equation}
The variational formulation of \eqref{eq:pde} reads as: find $u \in V$ such that
\begin{equation} \label{eq: formal var problem}
    \prescript{}{V^*}{\langle \mathcal{L}u, v \rangle}_{V} = \prescript{}{V^*}{\langle f, v \rangle}_{V} \quad \forall v \in V.
\end{equation}
The solution $u$ is formally expressed as $u = \mathcal{G} f$, with the Green's operator $\mathcal{G}: V^* \to V$, that is $\mathcal{G} = \mathcal{L}^{-1}$.

\noindent This formal problem may also be formulated in a variational multiscale setting. In order to do so, one introduces $\Bar{V}$ as a closed subspace of $V$, along with a linear projector $\mathcal{P}: V \rightarrow \Bar{V}$ with $\mathcal{P}^2 = \mathcal{P}$ and $\text{Range}(\mathcal{P}) = \Bar{V}$. $\mathcal{P}$ is assumed continuous in $V$. For the additional formal details, such as the required inf-sup conditions, see \cite{hughes_sangalli_2007}. Then, defining $V{'} = \text{Ker}(\mathcal{P})$, which is a closed subspace of $V$, combined with the continuity of $\mathcal{P}$ gives
\begin{equation}\label{eq: splitting V}
    V = \Bar{V} \oplus V{'}.
\end{equation}
This split implies that any $v \in V$ can be uniquely written as $v = \Bar{v} + v{'}$ with $\Bar{v} \in \Bar{V}$ and $v{'} \in V{'}$. 
In the VMS approach, $\Bar{V}$ represents the resolvable coarse scales, while $V{'}$ represents the unresolved fine scales. When using a VMS approach, one obtains $\Bar{u} = \mathcal{P}u$ as the solution for \eqref{eq:pde}, i.e. the numerical error equals the interpolation error and not a multiple of the interpolation error which is usually the case for optimal Finite Element methods. Using \eqref{eq: splitting V}, the variational formulation from \eqref{eq: formal var problem} can be split into
\begin{subequations}
    \begin{empheq}[left = \empheqlbrace\,]{align}
        \prescript{}{V^*}{\langle \mathcal{L}\Bar{u}, \Bar{v} \rangle}_{V} + \prescript{}{V^*}{\langle \mathcal{L}u{'}, \Bar{v} \rangle}_{V} &= \prescript{}{V^*}{\langle f, \Bar{v} \rangle}_{V} \quad \forall \Bar{v} \in \Bar{V}, \label{eq: cont var coarse}  \\
        \prescript{}{V^*}{\langle \mathcal{L}\Bar{u}, v{'} \rangle}_{V} + \prescript{}{V^*}{\langle \mathcal{L}u{'}, v{'} \rangle}_{V} &= \prescript{}{V^*}{\langle f, v{'} \rangle}_{V} \quad \forall v{'} \in V{'}.  \label{eq: cont var fine}
    \end{empheq}
\end{subequations}
These equations are assumed to be well-posed for $\Bar{u} \in \Bar{V}$ given $u{'}$ and $f$, and for $u{'} \in V{'}$ given $\Bar{u}$ and $f$, respectively. 

As for the non-split formulation in \eqref{eq: formal var problem}, a Green's operator can also be associated with the split formulation. In particular, the Greens' operator associated with the fine-scale equation in \eqref{eq: cont var fine} is the \emph{Fine-Scale Greens' operator} $\mathcal{G}{'}: V^* \to V{'}$. This gives $u{'}$ from the coarse-scale residual $\mathscr{R}\bar{u} = f - \mathcal{L}\Bar{u}$,
\begin{equation}
    u{'} = \mathcal{G}{'}(f - \mathcal{L}\Bar{u}) = \mathcal{G}{'}\mathscr{R}\bar{u}.
    \label{eq:uprine_def}
\end{equation}
With $\mathcal{G}{'}$, one can eliminate $u{'}$ from \eqref{eq: cont var coarse} and obtain the VMS formulation for $\Bar{u}$ as
\begin{equation} \label{eq: cont var full coarse}
    \prescript{}{V^*}{\langle \mathcal{L}\Bar{u}, \Bar{v} \rangle}_{V} - \prescript{}{V^*}{\langle \mathcal{L}\mathcal{G}{'}\mathcal{L}\Bar{u}, \Bar{v} \rangle}_{V} = \prescript{}{V^*}{\langle f, \Bar{v} \rangle}_{V} - \prescript{}{V^*}{\langle \mathcal{L}\mathcal{G}{'}f, \Bar{v} \rangle}_{V} \quad \forall \Bar{v} \in \Bar{V}.
\end{equation}
As a consequence of \eqref{eq: splitting V}, \eqref{eq: cont var full coarse} admits the unique solution $\Bar{u} = \mathcal{P}u$.

\noindent The above formulation can be used to obtain $\Bar{u}$, but requires obtaining the Fine-Scale Greens' operator $\mathcal{G}{'}$. It turns out (see \cite[\S 2.3]{hughes_sangalli_2007}) that this operator can be expressed in terms of the full Greens' operator, the projection, and the adjoint projection $\mathcal{P}^T : \Bar{V}^* \to V^*$ which is defined as
\begin{equation}
    \prescript{}{V^*}{\langle \mathcal{P}^T\Bar{\mu}, v \rangle}_{V} = \prescript{}{\Bar{V}^*}{\langle \Bar{\mu}, \mathcal{P}v \rangle}_{\Bar{V}} \quad \forall v \in V, \Bar{\mu} \in \Bar{V}^*,
\end{equation}
where $\Bar{V}^*$ is the dual space of $\Bar{V}$ and $\prescript{}{\Bar{V}^*}{\langle \cdot, \cdot \rangle}_{\Bar{V}}$ the pairing between them.

The Fine-Scale Greens' operator is then computed by considering a constrained version of \eqref{eq:uprine_def} as follows
\begin{gather}
    \mathcal{L} u' + \mathcal{P}^T \bar{\lambda} = \mathscr{R} \bar{u} \label{eq:fine_scale_con0}\\
    \mathcal{P} u' = 0, \label{eq:fine_scale_con1}
\end{gather}
where the newly introduced $\bar{\lambda} \in \bar{V}^*$ constrains $u{'}$ to live in the $V{'}$ space. Given the invertibility of $\mathcal{L}$, we may write
\begin{gather}
    u' = \mathcal{G}(\mathscr{R} \bar{u} - \mathcal{P}^T \bar{\lambda}),
\end{gather}
which when substituted into \eqref{eq:fine_scale_con1} and rearranged for $\bar{\lambda}$ gives
\begin{gather}
    \bar{\lambda} = \left(\mathcal{P} \mathcal{G} \mathcal{P}^T\right)^{-1} \mathcal{P} \mathcal{G} \mathscr{R} \bar{u}.
\end{gather}
The fine-scale operator $\mathcal{G}{'}$ can now be expressed in terms of $\mathcal{G}$ and $\mathcal{P}$ as
\begin{equation}\label{eq: g prime}
    \mathcal{G}{'} = \mathcal{G} - \mathcal{GP}^T (\mathcal{PGP}^T)^{-1}\mathcal{PG},
\end{equation}
with the properties
\begin{equation}
    \mathcal{G}{'}\mathcal{P}^T = 0, \quad \text{and} \quad \mathcal{PG}{'} = 0.
    \label{eq:g_p_prop}
\end{equation}
See \cite[p.542]{hughes_sangalli_2007} for the proof. \eqref{eq: g prime} can be used in \eqref{eq: cont var full coarse} to obtain 
\begin{equation}
    \prescript{}{V^*}{\langle \mathcal{L}\Bar{u}, \Bar{v} \rangle}_{V} - \prescript{}{V^*}{\langle \mathcal{L}\mathcal{G}{'}\mathcal{L}\Bar{u}, \Bar{v} \rangle}_{V} = \prescript{}{\Bar{V}^*}{\langle (\mathcal{PGP}^T)^{-1}\Bar{u}, \Bar{v} \rangle}_{\Bar{V}},
\end{equation}
which, since $(\mathcal{PGP}^T)^{-1}$ is invertible, confirms that \eqref{eq: cont var full coarse} is well-posed.
\subsection{Finite-dimensional case}
In the interest of practical implementations, $\bar{V}$ is a finite-dimensional subspace of $V$ and we seek a set of functionals $\boldsymbol{\mu} \in \Bar{V}^*$ which act as the basis for the projection $\mathcal{P}$. For an $N$ dimensional subspace $\bar{V}$ we seek $N$ functionals which satisfy
\begin{gather}
    \prescript{}{V^*}{\langle \mu_i, v \rangle}_{V} = 0 \quad \forall v \in V' \quad \forall i = 0, 1, \hdots, N - 1,
    \label{eq:mu_def}
\end{gather}
which in turn implies that $\mathcal{P}v = 0 \:\: \forall v \in V'$. These $\mu$'s form a basis for the range of $\mathcal{P}^T$ and we have 
\begin{gather}
    \mathcal{G}' \mu_i = 0 \quad \forall i = 0, 1, \hdots N - 1 \label{eq:Gprime_mu}\\
    \prescript{}{V^*}{\langle \mu_i, \mathcal{G}' \nu \rangle}_{V} = 0 \quad \forall \nu \in V^* \:\: \forall i = 0, 1, \hdots, N - 1, \label{eq:mu_Gprime_mu}
\end{gather}
following \eqref{eq:g_p_prop}. Furthermore, the $\mu$'s may be represented as a vector $\boldsymbol{\mu} \in (\Bar{V}^*)^N$
\begin{gather}
    \boldsymbol{\mu} = \left[\begin{array}{c}
         \mu_0 \\
         \mu_1 \\
         \vdots \\
         \mu_{N - 1}
    \end{array} \right]
    \quad \quad \text{and} \quad \quad \boldsymbol{\mu}^T = \left[\mu_0, \mu_1, \hdots, \mu_{N - 1} \right],
\end{gather}
through which $\mathcal{G}\boldsymbol{\mu}^T \in V^N$, $\boldsymbol{\mu} \mathcal{G} \boldsymbol{\mu}^T \in \mathbb{R}^{N \times N}$ and $\boldsymbol{\mu} \mathcal{G}$ can be written as 
\begin{gather}
    \mathcal{G}\boldsymbol{\mu}^T = \left[\mathcal{G} \mu_0, \mathcal{G} \mu_1, \hdots, \mathcal{G} \mu_{N - 1} \right] \\
    \boldsymbol{\mu} \mathcal{G} \boldsymbol{\mu}^T = \left[\begin{array}{ccc}
        \prescript{}{V^*}{\langle \mu_0, \mathcal{G} \mu_0 \rangle}_{V} & \hdots & \prescript{}{V^*}{\langle \mu_0, \mathcal{G} \mu_{N - 1} \rangle}_{V} \\
        \vdots & \ddots & \vdots \\
        \prescript{}{V^*}{\langle \mu_{N - 1}, \mathcal{G} \mu_0 \rangle}_{V} & \hdots & \prescript{}{V^*}{\langle \mu_{N - 1}, \mathcal{G} \mu_{N - 1} \rangle}_{V}
    \end{array} \right] \\
    \boldsymbol{\mu} \mathcal{G}(\nu) = \left[\begin{array}{c}
         \prescript{}{V^*}{\langle \mu_0, \mathcal{G} \nu \rangle}_{V} \\
         \vdots \\
         \prescript{}{V^*}{\langle \mu_{N - 1}, \mathcal{G} \nu \rangle}_{V}
    \end{array} \right] \quad \forall \nu \in V^*.
\end{gather}
The Fine-Scale Greens' function from \eqref{eq: g prime} may now be rewritten as
\begin{equation}
    \mathcal{G}' = \mathcal{G} - \mathcal{G}\boldsymbol{\mu}^T \left[\boldsymbol{\mu} \mathcal{G} \boldsymbol{\mu}^T \right]^{-1} \boldsymbol{\mu} \mathcal{G}.
    \label{eq:g_prime_mu}
\end{equation}

In theory, we are free to choose any set of $\boldsymbol{\mu}$ which satisfy \eqref{eq:mu_def}, \eqref{eq:Gprime_mu}, and \eqref{eq:mu_Gprime_mu}. In \cite{hughes_sangalli_2007} several definitions for $\boldsymbol{\mu}$ are proposed for obtaining the Fine-Scale Greens' function for the $H_0^1$ and $L^2$ projections of specific example problems. However, their generalisation for arbitrary projections or high-order spectral elements is not considered. This is precisely the knowledge gap that the present work attempts to fill. As such, we turn our attention to constructing a general definition for $\boldsymbol{\mu}$. Naturally, the choice of $\boldsymbol{\mu}$ depends on the projection as stated in \eqref{eq:mu_def}. Ideally, we want a generalised approach for defining the $\boldsymbol{\mu}$ corresponding to an arbitrary projector which then consequently yields a unique Fine-Scale Greens' function associated with that projection. To achieve this, we propose the use of dual basis functions. 

%% file: sections/04-dual_basis.tex
\section{Primal and Dual bases}
\label{sec:dual_basis}

The dual basis functions considered here were initially created in the context of the Mimetic Spectral Element Method (MSEM) \cite{jain_zhang_palha_gerritsma_2021}. These are dual functions in the $L^2$-sense. In order to understand their link to projections and the use of dual basis functions as the $\boldsymbol{\mu}$'s, we first introduce the primal basis functions and their corresponding $L^2$ duals. For a more detailed exposition and the extensions to multiple dimensions, see \cite{jain_zhang_palha_gerritsma_2021}.

\subsection{Primal basis}\label{subsec:primal_basis}
\begin{figure}[H]
\begin{subfigure}{0.49\linewidth}
    \centering
    \includegraphics[width = \linewidth]{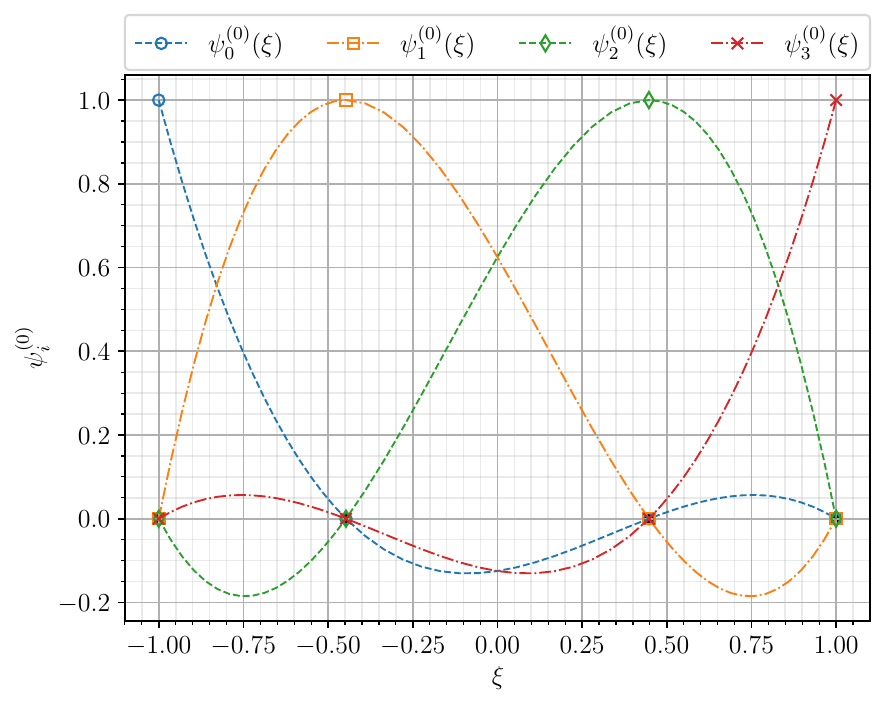}
    \caption{Nodal Lagrange polynomials of degree $p = 3$}
    \label{fig:nodal_poly}
\end{subfigure}
\begin{subfigure}{0.49\linewidth}
    \centering
    \includegraphics[width = \linewidth]{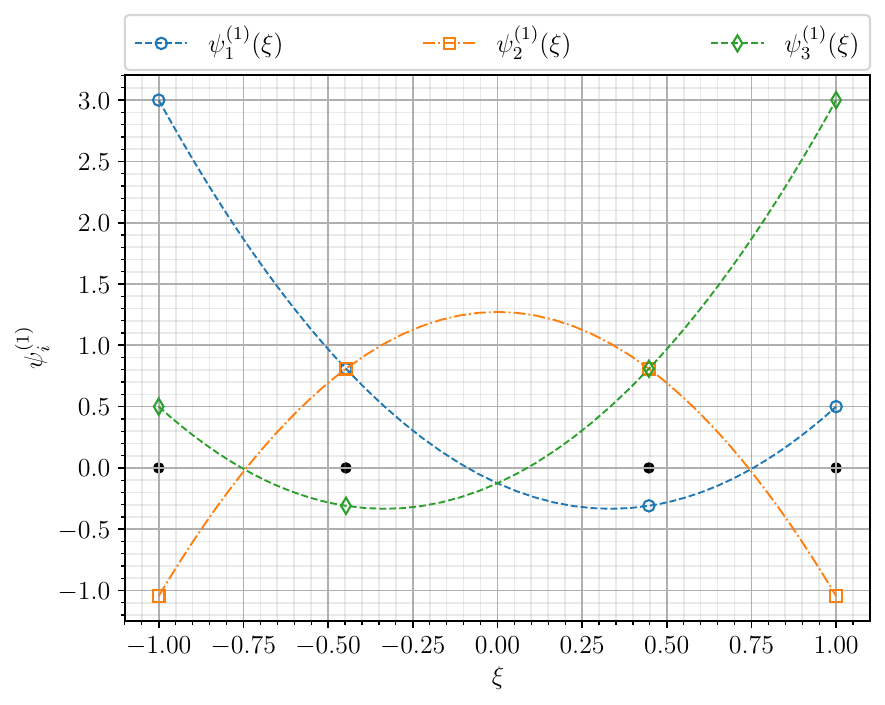}
    \caption{Edge polynomials of degree $p = 3$}
    \label{fig:edge_poly}
\end{subfigure}
\caption{Nodal Lagrange and corresponding Edge polynomials constructed over 1D mesh with 4 GLL points}
\end{figure}

A primal basis for nodal functions is constructed based on nodal degrees of freedom. Given such degrees of freedom, denoted by $\mathcal{N}^0_j$, the primal basis functions, denoted by $\psi^{(0)}_i$, should satisfy the following property
\begin{equation}
    \mathcal{N}^0_j(\psi_i^{(0)}) = \delta_{ij}.
    \label{eq:nodal_prop}
\end{equation}
We consider the mimetic spectral element method which uses the Gauss-Legendre-Lobatto (GLL) nodes for nodal degrees of freedom. To define the GLL-nodes, consider the interval $K = [-1, 1] \subset \mathbb{R}$ and let $\xi_i \in K, i = 0, \dots, p$ be the roots of the polynomial $\left(1 - \xi^2\right)L'_p(\xi)$, where $L_p(\xi)$ is the Legendre polynomial of degree $p$ and $L'_p(\xi)$ its derivative. The GLL-nodes are then used to define basis functions that satisfy \eqref{eq:nodal_prop} with $\mathcal{N}_j^{(0)}(f):=f(\xi_j)$
\begin{equation}
    \psi^{(0)}_i(\xi) := \displaystyle\prod_{\shortstack{$\scriptstyle j = 0$\\$\scriptstyle j \neq i$}}^{p} \frac{\xi - \xi_j}{\xi_i - \xi_j}, \quad \quad i = 0, 1, \hdots, p\;.
    \label{eq:lagrange_p}
\end{equation}
The basis functions $\psi_i^{(0)}(\xi)$, $i=0,\ldots,p$ span the finite dimensional subspace $\bar{V}$.
An example of these basis functions for $p=3$ is shown in \Cref{fig:nodal_poly}.

\bigskip
\noindent Next to the basis functions for nodal functions, the MSEM also defines basis functions based on edge degrees of freedom. These are histopolant basis functions based on integral degrees of freedom, defined as
\begin{equation}
    \mathcal{N}_j^1 (g) = \int_{\xi_{j-1}}^{\xi_j} g(\xi) \,d\xi, \quad i = 1, \dots, p,
\end{equation}
where $\mathcal{N}_{j}^{1}$ is a linear functional which extracts integral values from any function $g(\xi)$

As with the nodal basis, in order to define an edge basis, one looks for basis functions, now denoted $\psi^{(1)}$, that satisfy
\begin{equation}
    \mathcal{N}^1_j(\psi^{(1)}_i) = \int_{\xi_{j-1}}^{\xi_j} \psi^{(1)}_i(\xi) \,d\xi = \delta_{ij}.
    \label{eq:edge_prop}
\end{equation}
The edge basis functions on the GLL-grid can be expressed in terms of the nodal basis as \cite{gerritsma_2010}
\begin{equation}
    \psi^{(1)}_i(\xi) := -\sum_{k = 0}^{i - 1} \frac{\mathrm{d} \psi^{(0)}_k}{\mathrm{d}{\xi}}, \quad \quad i = 1, 2, \hdots, p.
\end{equation}
An example of these basis functions for $p=3$ is shown in \Cref{fig:edge_poly}. Note that the space spanned by these edge functions is the space of polynomials of degree $(p-1)$. 
\begin{remark}
    All subsequent references to edge polynomials denoted with degree $p$ refer to a polynomial basis of degree $(p - 1)$.
    \label{rem:edge_p}
\end{remark}

\subsection{Dual basis} \label{sec:dualSpace}
Based on the primal basis introduced in the previous section, one can derive an associated dual basis. Before doing so, we first review the definition of a dual space.
\begin{definition}[\cite{frankel}]
    \textit{A \emph{linear functional} $\alpha$ defined on a vector space $\mathcal{E}$ is a real-valued linear function that is a mapping $\alpha : \mathcal{E} \xrightarrow[]{} \mathbb{R}$}
    \begin{equation}
        \alpha(a\boldsymbol{u} + b\boldsymbol{v}) = a \alpha(\boldsymbol{u}) + b \alpha(\boldsymbol{v}), \quad \quad \forall \: \boldsymbol{u}, \boldsymbol{v} \in \mathcal{E} \: \text{ and } \: \forall \: a, b \in \mathbb{R}.
    \end{equation}
    \textit{If $\mathcal{E}$ is a finite-dimensional vector space with $n$ dimensions which has $\boldsymbol{e}_1, \boldsymbol{e}_2, \hdots, \boldsymbol{e}_{n}$ as its basis, any arbitrary vector $\boldsymbol{v} \in \mathcal{E}$ can be expanded as a linear combination of the basis vectors; $\boldsymbol{v} = \displaystyle \sum_{j} v^j \boldsymbol{e}_j$. Subsequently, applying the functional $\alpha$ to $\boldsymbol{v}$ yields:}
    \begin{equation}
        \alpha(\boldsymbol{v}) = \alpha \left(\sum_{j} v^j \boldsymbol{e}_j\right) = \sum_j v^j \alpha(\boldsymbol{e}_j).
        \label{eq:linFunctional}
    \end{equation}
\end{definition}
Note that the linear functional $\alpha$ does \emph{not} live in the vector space $\mathcal{E}$. 
Instead, the linear functional lives in a different space defined as follows.
\begin{definition}[\cite{frankel}]
    \textit{The collection of all the linear functionals on the vector space $\mathcal{E}$ form a separate vector space referred to as the \emph{dual space} denoted by $\mathcal{E^*}$. If $\mathcal{E}$ is n-dimensional, then so is $\mathcal{E^*}$, which enables one to define the \emph{dual basis} $\sigma^1, \sigma^2, \hdots, \sigma^n$. 
    The elements in the dual space are referred to as \emph{covectors}. Furthermore, the canonical dual basis satisfies the property:}
    \begin{equation}
        \sigma^i(\boldsymbol{e}_j) = \delta^{i}_{j},
        \label{eq:dual_delta}
    \end{equation}
    \textit{which, along with the linearity of the $\sigma^i$, can be used to show that applying $\sigma^i$ to a vector $\boldsymbol{v} \in \mathcal{E}$ simply extracts the $i^{th}$ component of $\boldsymbol{v}$}
    \begin{equation}
        \sigma^i\left(\sum_j \boldsymbol{e}_j v^j \right) = \sum_j \sigma^i (\boldsymbol{e}_j) v^j = \sum_j \delta^{i}_j v^j = v^i.
    \end{equation}
\end{definition}

The above definition can be used to construct a discrete dual space based on the primal spaces of \Cref{subsec:primal_basis}. For the MSEM, the $L^2$ inner product was used to create such a space \cite{jain_zhang_palha_gerritsma_2021}, which will be briefly explained in the next section.


\subsection{\texorpdfstring{$L^2$}{L2} dual basis}
To construct the dual bases, we represent $p^h$ with respect to a nodal basis as
\begin{equation}
    p^h (\xi) = \sum_{i=0}^p \mathcal{N}_i^0 (p^h) \psi_i^{(0)}(\xi).
\end{equation}
In this equation, $\mathcal{N}_i^0 (p^h)$ are the nodal degrees of freedom, and $\psi^{(0)}_i(\xi)$ are the associated basis functions (see \eqref{eq:lagrange_p}). For clarity, the following shorthand notation will be used \cite{jain_zhang_palha_gerritsma_2021},
\begin{equation}\label{eq:phExpansion}
    p^h(\xi) = \boldsymbol{\psi}^{(0)}(\xi) \mathcal{N}^0 (p^h),
\end{equation}
where
\begin{equation}
    \boldsymbol{\psi}^{(0)} (\xi) = 
    \left [
    \psi^{(0)}_0 (\xi) \;\;  \psi^{(0)}_1 (\xi) \;\;  \cdots \;\; \psi^{(0)}_{p-1} (\xi) \;\; \psi^{(0)}_p (\xi)
    \right ]
    \quad
    \text{and}
    \quad
    \mathcal{N}^0 (p^h) = 
    \left [ \begin{array}{c}
    \mathcal{N}_0^0 (p^h) \\[6pt] 
    \mathcal{N}_1^0 (p^h) \\ 
    \vdots \\ 
    \mathcal{N}_{p-1}^0 (p^h) \\[6pt] 
    \mathcal{N}_p^0 (p^h)
    \end{array} \right ] .
\end{equation}

\begin{figure}[ht]
\begin{subfigure}{0.49\linewidth}
    \centering
    \includegraphics[width = \linewidth]{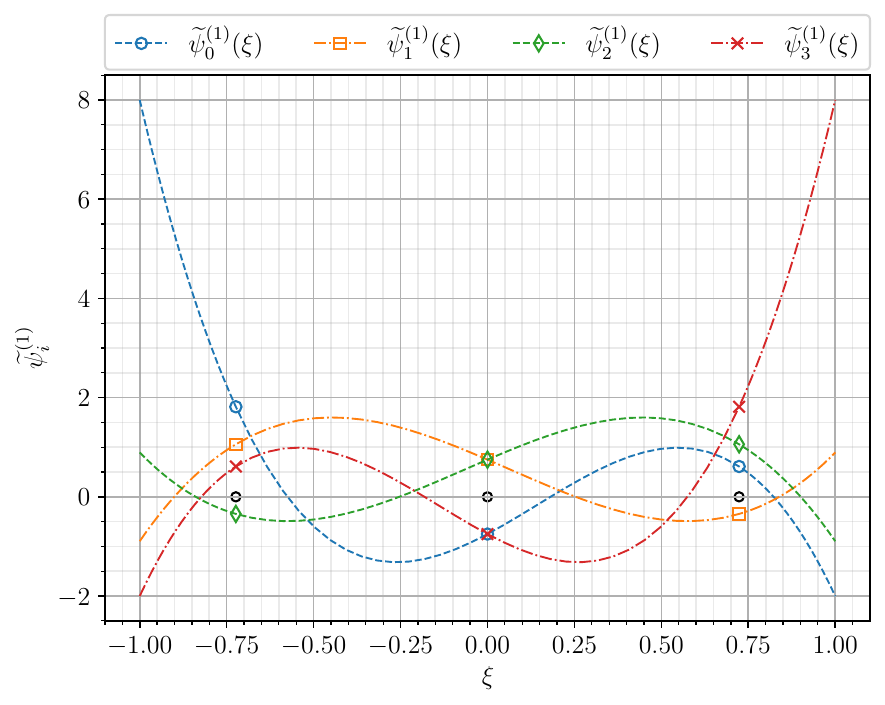}
    \caption{Dual edge polynomials of degree $p = 3$}
    \label{fig:dualEdge_poly}
\end{subfigure}
\begin{subfigure}{0.49\linewidth}
    \centering
    \includegraphics[width = \linewidth]{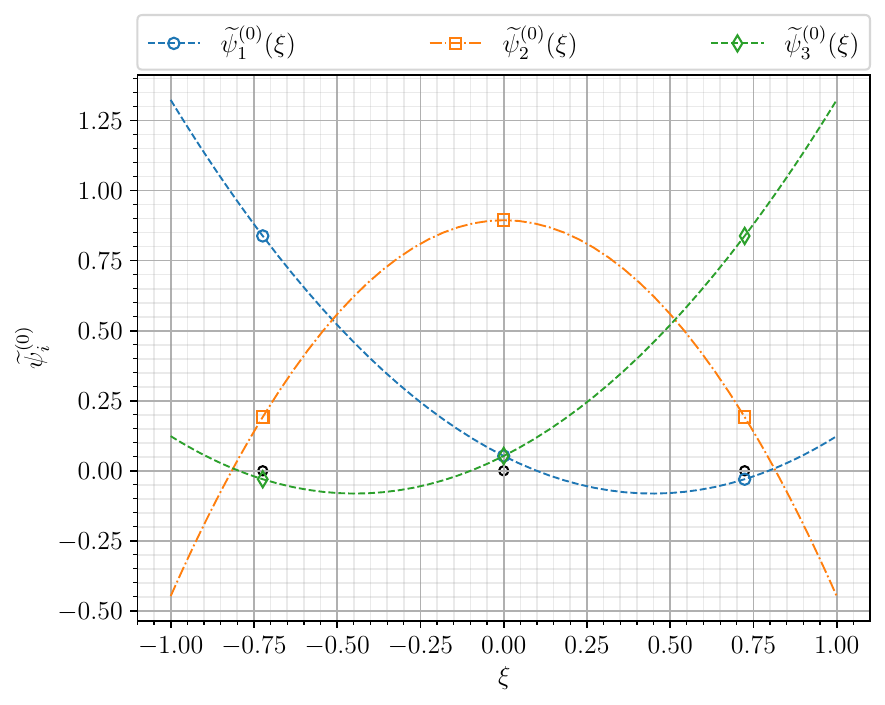}
    \caption{Dual nodal polynomials of degree $p = 3$}
    \label{fig:dualNodal_poly}
\end{subfigure}
\caption{Dual nodal and edge polynomials constructed over 1D mesh with 4 GLL points}
\end{figure}

Now, let $p^h$ and $q^h$ be expanded as in \eqref{eq:phExpansion}. The $L^2$-inner product between the two is then given by
\begin{equation}
    \left( p^h, q^h \right)_{L^2(K)} = \int_K p^h(\xi) q^h(\xi)\,\mathrm{d}K = \mathcal{N}^0(p^h)^T \MassM{0} \mathcal{N}^0(q^h),
\end{equation}
where $\MassM{0}$ denotes the mass matrix built using the nodal basis functions as 
\begin{equation}
    \MassM{0} = \int_K \boldsymbol{\psi}^{(0)}(\xi)^T \boldsymbol{\psi}^{(0)}(\xi)\,\mathrm{d}K.
\end{equation}
This mass matrix, and its inverse, can be used to define the dual degrees of freedom and the associated dual basis. The dual edge degrees of freedom (dual to the primal nodes) are then defined as
\begin{equation}
    \widetilde{\mathcal{N}}^1(p^h) = \MassM{0} \mathcal{N}^0(p^h).
\end{equation}
Note that these are the degrees of freedom for the dual line elements. The dual basis functions are then
\begin{equation}
    \widetilde{\boldsymbol{\psi}}^{(1)}(\xi) = \boldsymbol{\psi}^{(0)}(\xi) \MassMinv{0}.
\end{equation}
See \cite{jain_zhang_palha_gerritsma_2021} for the details and proofs. The plots of these dual edge polynomials are shown in \Cref{fig:dualEdge_poly}.

Similarly, the basis dual to the primal edge polynomials can also be constructed. Defining $p^h$ and $q^h$ as 
\begin{equation}
    p^h(\xi) = \boldsymbol{\psi}^{(1)}(\xi) \mathcal{N}^1 (p^h), \quad \quad q^h(\xi) = \boldsymbol{\psi}^{(1)}(\xi) \mathcal{N}^1 (q^h)
\end{equation}
with 
\begin{equation}
    \boldsymbol{\psi}^{(1)} (\xi) = 
    \left [
    \psi^{(1)}_1 (\xi) \;\; \psi^{(1)}_2 (\xi) \;\; \cdots \;\; \psi^{(1)}_{p-1} (\xi) \;\; \psi^{(1)}_p (\xi)
    \right ]
    \quad
    \text{and}
    \quad
    \mathcal{N}^1 (p^h) = 
    \left [ \begin{array}{c}
    \mathcal{N}_1^1 (p^h) \\[6pt] 
    \mathcal{N}_2^1 (p^h) \\ 
    \vdots \\ 
    \mathcal{N}_{p-1}^1 (p^h) \\[6pt] 
    \mathcal{N}_p^1 (p^h)
    \end{array} \right ] .
\end{equation}
The $L^2$-inner product is now given by
\begin{equation}
    \left( p^h, q^h \right)_{L^2(K)} = \int_K p^h(\xi) q^h(\xi)\,dK = \mathcal{N}^1(p^h)^T \mathbb{M}^1 \mathcal{N}^1(q^h),
\end{equation}
where $\MassM{1}$ denotes the mass matrix built using the edge basis functions as 
\begin{equation}
    \MassM{1} = \int_K \boldsymbol{\psi}^{(1)}(\xi)^T \boldsymbol{\psi}^{(1)}(\xi)\,dK.
\end{equation}
Constructing the dual degrees of freedom and dual basis follows the exact same procedure as in the nodal case. Thus, the dual nodal degrees of freedom (dual to the primal edges) are defined as
\begin{equation}
    \widetilde{\mathcal{N}}^0(p^h) = \MassM{1} \mathcal{N}^1(p^h).
    \label{eq:dual_0_dof}
\end{equation}
Note that these are the degrees of freedom for the dual nodal elements. The dual basis functions are then
\begin{equation}
    \widetilde{\boldsymbol{\psi}}^{(0)}(\xi) = \boldsymbol{\psi}^{(1)}(\xi) \MassMinv{1}.
    \label{eq:dual_0_basis}
\end{equation}
See again \cite{jain_zhang_palha_gerritsma_2021} for the details and proofs. The plots of these dual polynomials are shown in \Cref{fig:dualNodal_poly}.

The primal and dual bases discussed thus far have all been defined in a single (reference) domain $\xi \in K = [-1, 1]$. For the generalisation to a multi-element setting over an arbitrary physical domain $x \in \Omega = [a, b]$ divided into $\mathrm{N}$ elements with $a, b \in \mathbb{R}$ and $a = x_0 < x_1 < \hdots < x_{\mathrm{N} - 1} < x_{\mathrm{N}} = b$, we introduce the following linear map $\Phi$ between $K$ and $[x_{n - 1}, x_n]$
\begin{equation}
    x = \Phi(\xi) = \frac{1}{2}(1 - \xi) x_{n - 1} + \frac{1}{2} (1 + \xi) x_n, \quad \quad -1 \leq \xi \leq 1.
\end{equation}
Subsequently, the global nodal basis functions are defined as
\begin{equation}
    \psi^{(0)}_i(x) = \psi^{(0)}_j(\xi) \circ \Phi^{-1}(x) \quad \text{with } i = j + (n - 1)p, 
\end{equation}
where $n = 1, \hdots, \mathrm{N}$ and $j = 0, \hdots p$, giving a total of $k = p + 1 + (\mathrm{N} - 1) p$ nodal functions. Similarly, the global edge basis functions are given by
\begin{equation}
    \psi^{(1)}_i(x) = \psi^{(1)}_j(\xi) \circ \frac{\Phi^{-1}(x)}{\boldsymbol{J}} \quad \text{with } i = j + (n - 1)p, 
\end{equation}
where $n = 1, \hdots, \mathrm{N}$ and $j = 1, \hdots p$, and $\boldsymbol{J} = \frac{\mathrm{d}x}{\mathrm{d}\xi}$ is the Jacobian, giving a total of $k = p + (\mathrm{N} - 1) p$ edge functions. Further details regarding the multi-element case may be found in \cite[\S 2.3]{jain_zhang_palha_gerritsma_2021}.



\subsection{Projections}
We define the projections as follows
\begin{definition}
    \textit{Consider a real-valued infinite-dimensional function space $V$ and a finite-dimensional function subspace $\bar{V}$. A \emph{projection} is the mapping $\mathcal{P} : V \xrightarrow[]{} \bar{V}$ from the infinite-dimensional space to the finite-dimensional one which minimises the difference between the infinite and finite-dimensional representations in a given norm.}
\end{definition}

\subsubsection{\texorpdfstring{$L^2$}{L2} projection}\label{sec:L2_projection}
\begin{figure}[H]
\begin{subfigure}{0.49\linewidth}
    \centering
    \includegraphics[width = \linewidth]{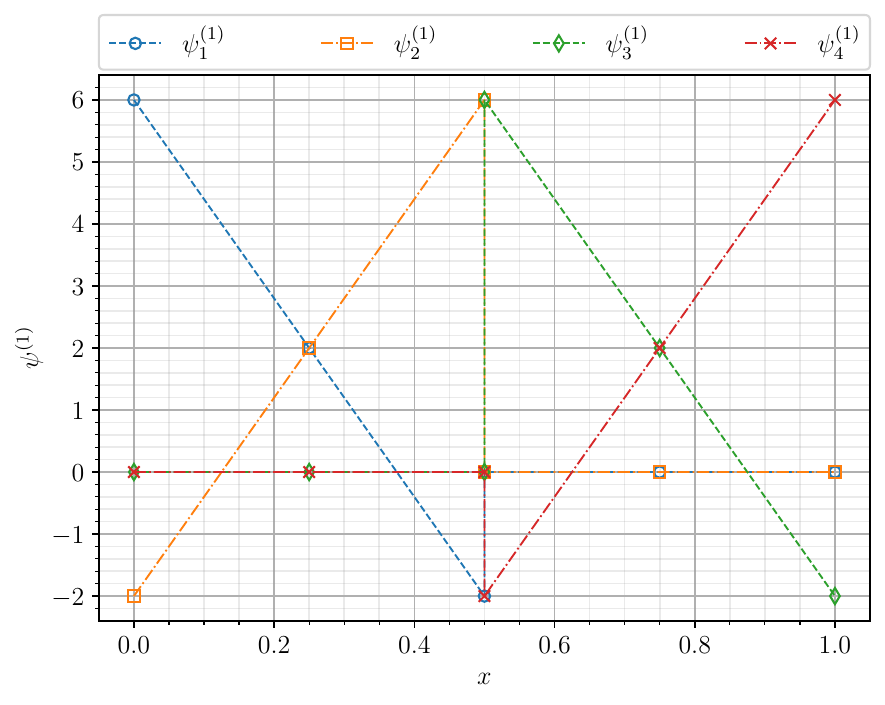}
    \caption{Edge basis functions}
\end{subfigure}
\begin{subfigure}{0.49\linewidth}
    \centering
    \includegraphics[width = \linewidth]{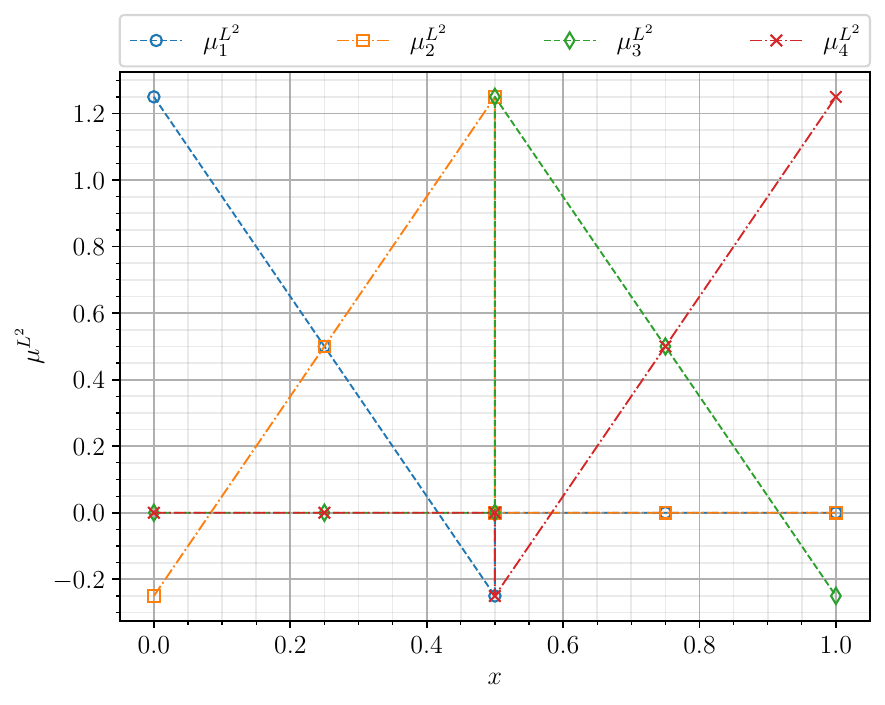}
    \caption{$L^2$ dual basis functions}
\end{subfigure}
\caption{Primal edge basis functions of polynomial degree $p = 2$ over 2 elements and their corresponding $L^2$ duals over $x \in \Omega \in [0, 1]$}
\label{fig:L2_dual_p1}
\end{figure}

\begin{figure}[H]
\begin{subfigure}{0.49\linewidth}
    \centering
    \includegraphics[width = \linewidth]{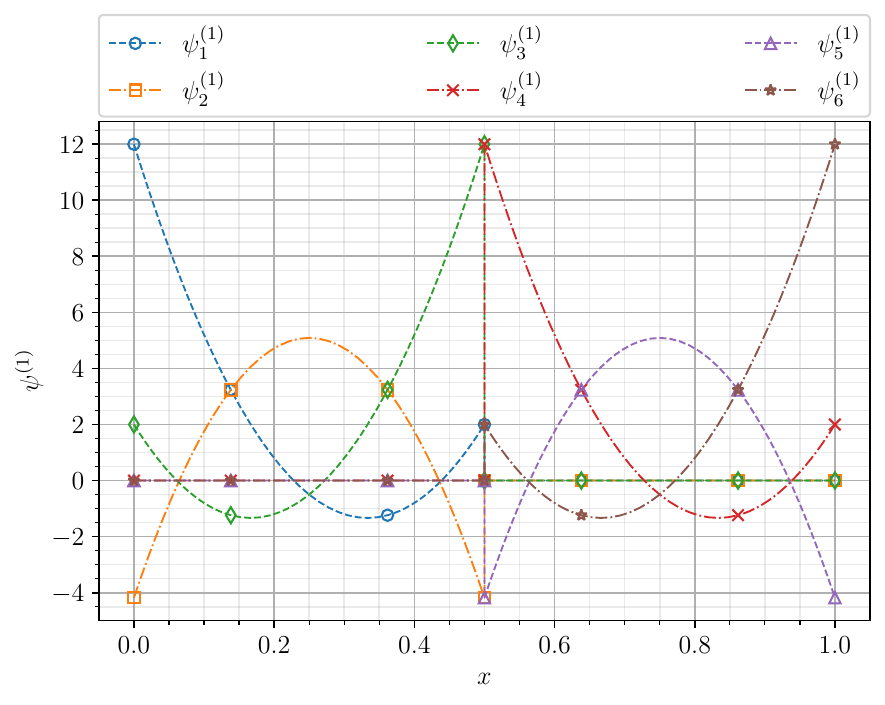}
    \caption{Edge basis functions}
\end{subfigure}
\begin{subfigure}{0.49\linewidth}
    \centering
    \includegraphics[width = \linewidth]{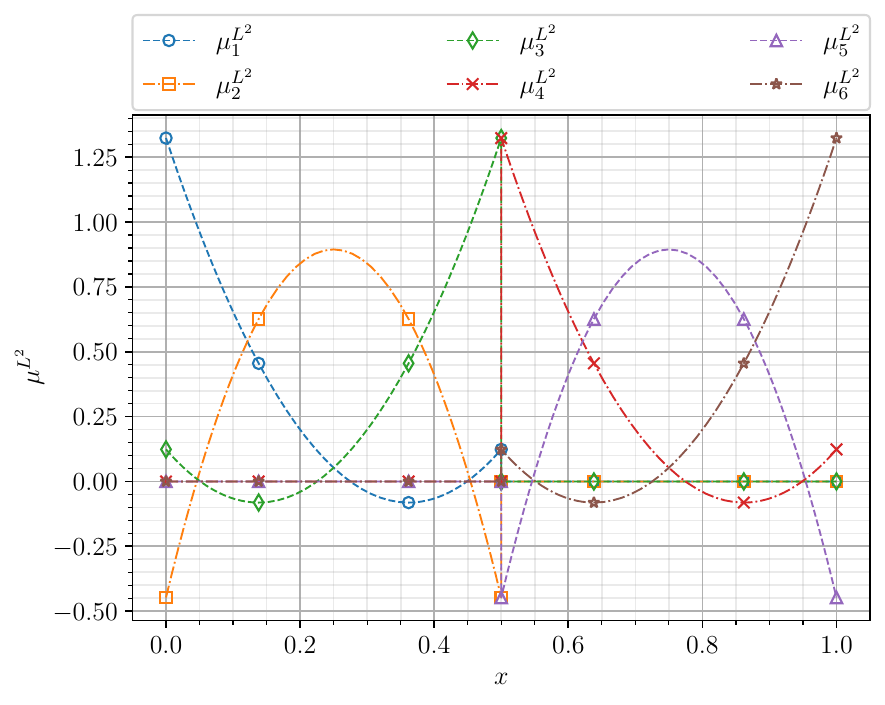}
    \caption{$L^2$ dual basis functions}
\end{subfigure}
\caption{Primal edge basis functions of polynomial degree $p = 3$ over 2 elements and their corresponding $L^2$ duals over $x \in \Omega \in [0, 1]$}
\label{fig:L2_dual_p2}
\end{figure}

Effectively, the projection maps continuous quantities (functions) into the discrete finite-dimensional subspace such that the best possible discrete representation in the considered norm is obtained. We will show that dual basis functions encode the projector. To demonstrate this fact, consider the $L^2$ projection $\mathcal{P}_{L^2}$ of a function $\phi \in L^2(\Omega)$.
The minimisation problem thus reads
\begin{gather}
    \mathcal{P}_{L^2} \phi = \argmin_{\bar{\phi} \in \bar{V}} \left\{\frac{1}{2} \left\lVert \bar{\phi} - \phi \right\rVert_{L^2}^2\right\}\;,
\end{gather}
which gives
\begin{gather}
    \left(v^h, \bar{\phi} \right)_{L^2(\Omega)} = \left(v^h, \phi \right)_{L^2(\Omega)}, \quad \forall v^h \in \bar{V}. \label{eq:L2_proj}
\end{gather}
Considering a 1D case and expanding $\bar{\phi}$ in terms of primal edge polynomials and $v_h$ in terms of the corresponding dual basis given by \eqref{eq:dual_0_dof} and \eqref{eq:dual_0_basis}
\begin{gather}
    \bar{\phi} := \boldsymbol{\psi}^{(1)} \mathcal{N}^1(\bar{\phi}) \label{eq:expansion_projection}\\
    v^h := \widetilde{\boldsymbol{\psi}}^{(0)} \widetilde{\mathcal{N}}^{0}(v^h),
\end{gather}
we obtain after substitution of these expansions into \eqref{eq:L2_proj},
\begin{gather}
    \widetilde{\mathcal{N}}^{0}(v^h)^T \left(\widetilde{\boldsymbol{\psi}}^{(0)}, \boldsymbol{\psi}^{(1)} \right)_{L^2(\Omega)} \mathcal{N}^1(\bar{\phi}) = \widetilde{\mathcal{N}}^{0}(v^h)^T \left(\widetilde{\boldsymbol{\psi}}^{(0)}, \phi \right)_{L^2(\Omega)}, \quad \forall \widetilde{\mathcal{N}}^{0}(v^h) \in \mathbb{R}^{k}.
\end{gather}
Using the fact that the primal and the dual bases are bi-orthogonal, i.e. $\left(\widetilde{\boldsymbol{\psi}}^{(0)}_i, \boldsymbol{\psi}^{(1)}_j \right)_{L^2(\Omega)}=\delta_{i,j}$, we have
\begin{gather}
    \mathcal{N}^1(\bar{\phi}) = \left(\boldsymbol{\mu}^{L^2}, \phi \right)_{L^2(\Omega)},
\end{gather}
where we introduced $\boldsymbol{\mu}^{L^2}=\widetilde{\boldsymbol{\psi}}^{(0)}$.
This indicates that the inner product between the $L^2$ dual basis and any $\phi \in L^2(\Omega)$ returns the expansion coefficients for the projection of $\phi$ onto $\bar{V}$. The projected solution is then obtained using the primal edge functions, \eqref{eq:expansion_projection}. The plots of these $L^2$ dual functions in a multi-element setting are shown in \Cref{fig:L2_dual_p1} and \Cref{fig:L2_dual_p2} for different polynomial degrees. 

\begin{remark}
    We choose the space of edge basis functions to expand the functions in $\bar{V}$ for the $L^2$ projection for the reason that the edge functions live strictly in the $L^2$ space. We may, however, also choose to expand $\bar{\phi}$ in nodal basis functions which would give us the nodal $L^2$ projection.
    \label{rem:L2_space}
\end{remark}




\subsubsection{\texorpdfstring{$H^1_0$}{H10} projection}\label{sec:H1_pojection}

The concept introduced need not be limited to $L^2$. This may also be extended to other projections giving rise to different sets of dual basis functions. For instance, considering the $H_0^1$ projection yields the $H_0^1$ dual basis functions. We consider the $H_0^1$ projection in a manner similar to that for the $L^2$ but with the projector $\mathcal{P}_{H_0^1}$.
For $\boldsymbol{\phi}\in H^1_0(\Omega)$ the minimisation problem reads
\begin{gather}
    \mathcal{P}_{H_0^1} \phi = \argmin_{\bar{\phi} \in \bar{V}} \left\{\frac{1}{2} \left\lVert \nabla \bar{\phi} - \nabla \phi \right\rVert_{L^2}^2\right\},
\end{gather}
which leads to
\begin{gather}
    \left(v^h, \bar{\phi}\right)_{H^1_0(\Omega)} = \left(v^h, \phi \right)_{H^1_0(\Omega)}, \quad \forall v^h \in \bar{V}. \label{eq:H01_opt}
\end{gather}
For the 1D case, we substitute
\begin{gather}
    \bar{\phi} := \boldsymbol{\psi}^{(0)} \mathcal{N}^0(\bar{\phi}) \label{eq:expansion_phi_bar}\\
    -\Delta_h v^h := \boldsymbol{\psi}^{(0)} \mathcal{N}^0(v^h) \quad \Longleftrightarrow \quad v^h = \left ( -\Delta^{-1}_h \boldsymbol{\psi}^{(0)}  \right ) \mathcal{N}^0(v^h),
\end{gather}
into \eqref{eq:H01_opt} where $\Delta^{-1}_h$ is the \textit{discrete} inverse Laplacian operator with homogeneous Dirichlet boundaries, and we get
\begin{gather}
    \mathcal{N}^0(v^h)^T \left(-\Delta^{-1}_h\boldsymbol{\psi}^{(0)} , \boldsymbol{\psi}^{(0)}\right)_{H^1_0(\Omega)} \mathcal{N}^0(\bar{\phi}) = \mathcal{N}^0(v^h)^T  \left(-\Delta^{-1}_h\boldsymbol{\psi}^{(0)} , \phi \right)_{H^1_0(\Omega)}, \quad \forall \mathcal{N}^0(v^h)  \in \mathbb{R}^{k} \;. 
\end{gather}
Using the fact that $\left(-\Delta^{-1}_h \boldsymbol{\psi}^{(0)}_i , \boldsymbol{\psi}^{(0)}_j \right)_{H^1_0(\Omega)}=\delta_{i,j}$ and the fact that this equation needs to hold for all $\mathcal{N}^0(v^h)  \in \mathbb{R}^{k}$ gives
\begin{equation}
    \mathcal{N}^0(\bar{\phi}) = \left(-\Delta^{-1}_h \boldsymbol{\psi}^{(0)} , \phi \right)_{H^1_0(\Omega)} \;.
\end{equation}
If we define $\boldsymbol{\mu}^{H^{1}_0}$ to be the solution of $-\Delta_h \boldsymbol{\mu}^{H^{1}_0} = \boldsymbol{\psi}^{(0)}$, this can be written as
\begin{equation}
    \mathcal{N}^0(\bar{\phi}) = \left(\boldsymbol{\mu}^{H^{1}_0}, \phi \right)_{H^1_0(\Omega)}. \label{eq:u_H01_proj_dual}
\end{equation}
Once again, when the expansion coefficients $\mathcal{N}^0(\bar{\phi})$ are obtained from \eqref{eq:u_H01_proj_dual}, and we use \eqref{eq:expansion_phi_bar} to produce the $H^1_0$ projected solution in $\bar{V}$. 

\begin{remark}
    The functionals $\boldsymbol{\mu}$ defined in Sections~\ref{sec:L2_projection} and \ref{sec:H1_pojection}, are explicitly computable, and play the role of the $\{ \mu_i \}_{i=1,\ldots,N}$ described in the paper of Hughes and Sangalli, \cite[(2.17)]{hughes_sangalli_2007}.
\end{remark}

\begin{figure}[H]
\begin{subfigure}{0.49\linewidth}
    \centering
    \includegraphics[width = \linewidth]{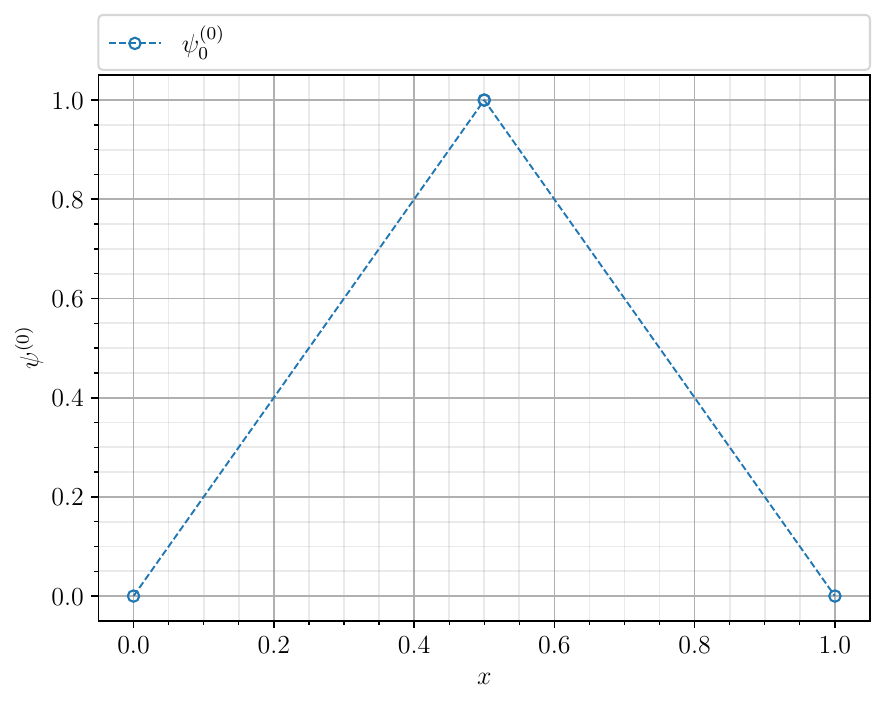}
    \caption{Nodal basis function}
\end{subfigure}
\begin{subfigure}{0.49\linewidth}
    \centering
    \includegraphics[width = \linewidth]{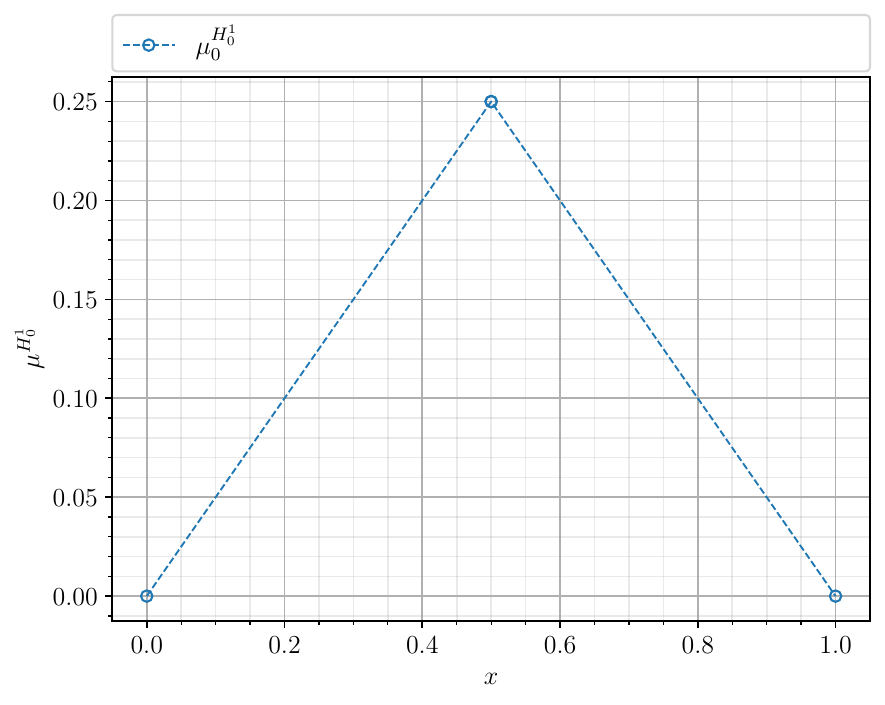}
    \caption{$H_0^1$ dual basis function}
\end{subfigure}
\caption{Primal nodal basis function of polynomial degree $p = 1$ over 2 elements and its corresponding $H_0^1$ dual over $x \in \Omega \in [0, 1]$}
\label{fig:H01_dual_p1}
\end{figure}

\begin{figure}[H]
\begin{subfigure}{0.49\linewidth}
    \centering
    \includegraphics[width = \linewidth]{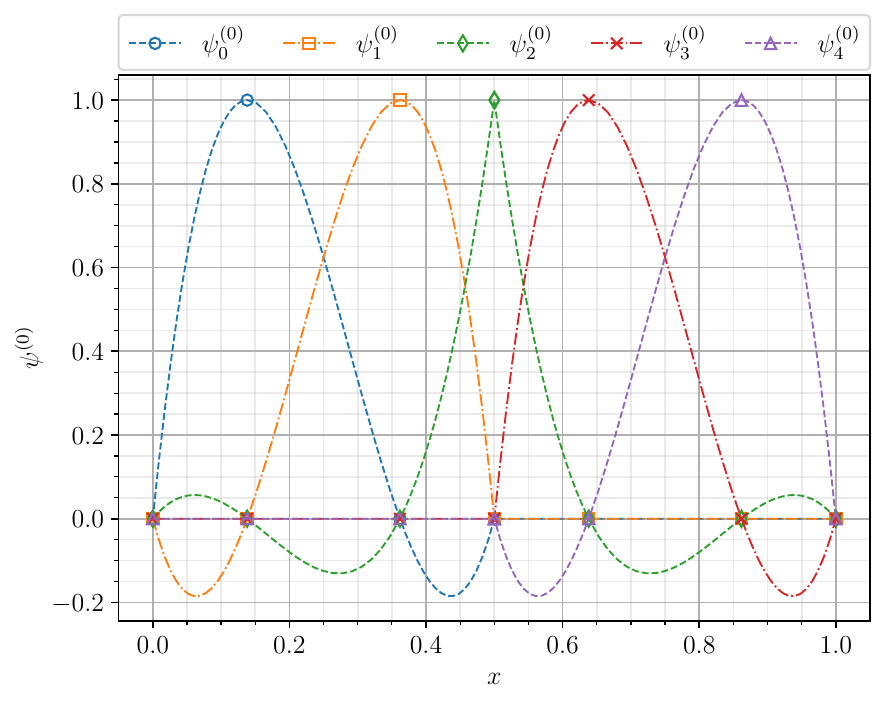}
    \caption{Nodal basis functions}
\end{subfigure}
\begin{subfigure}{0.49\linewidth}
    \centering
    \includegraphics[width = \linewidth]{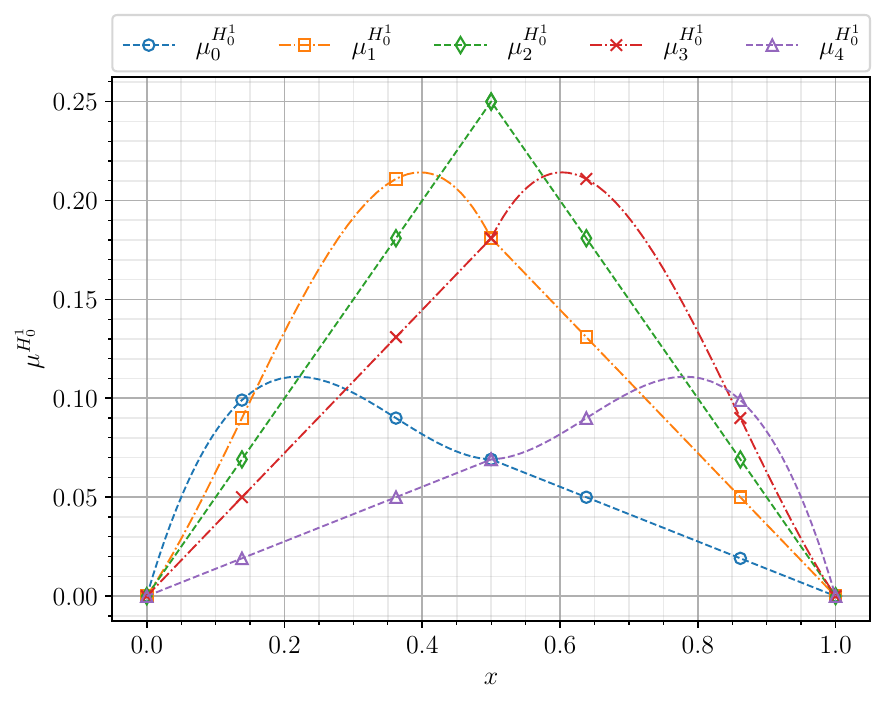}
    \caption{$H_0^1$ dual basis functions}
\end{subfigure}
\caption{Primal nodal basis functions of polynomial degree $p = 3$ over 2 elements and their corresponding $H_0^1$ duals over $x \in \Omega \in [0, 1]$}
\label{fig:H01_dual_p3}
\end{figure}

\noindent Plots of these dual bases are presented in \Cref{fig:H01_dual_p1} and \Cref{fig:H01_dual_p3}. It is evident that the $H_0^1$ dual functions are all positive functions and inherently global.

%% file: sections/05-greens_function.tex
\section{Fine-Scale Greens' function}
\label{sec:greens}



Having reviewed the formulation of the Fine-Scale Greens' function and proposing to define the $\boldsymbol{\mu}$'s as dual basis functions, we now proceed to construct the Fine-Scale Greens' function. In this section, we present the various computations involved in constructing the Fine-Scale Greens' function for the 1D Poisson equation associated with both the $H_0^1$ and $L^2$ projections. 

\subsection{Computing the Fine-Scale Greens' function for 1D Poisson}
\begin{figure}[H]
    \centering
    \includegraphics[width = 0.45\linewidth]{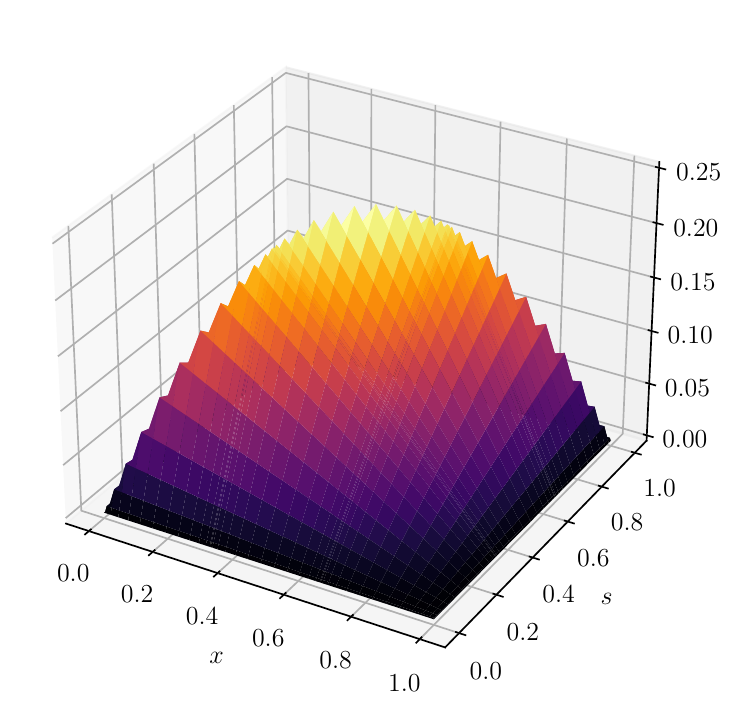}
    \caption{Greens' function for the 1D Poisson in $\Omega = [0, 1]$}
    \label{fig:poisson_greens}
\end{figure}

The 1D Poisson problem is stated as follows
\begin{gather}
    -\parddxdx{u}{x} = f, \quad x \in \Omega = [0, 1] \label{eq:poisson_eq}\\
    u(x) = 0, \quad x \in \partial \Omega. \label{eq:poisson_bc}
\end{gather}
The (global) Greens' function for this problem defined in $\Omega$ is given by
\begin{equation}
    g(x, s) = \begin{cases}
        (1 - s)x, \quad x \leq s \\
        s (1 - x), \quad x > s
    \end{cases},
\end{equation}
the plot of which is shown in \Cref{fig:poisson_greens}.
Using the Greens' function and the $\boldsymbol{\mu}$ defined by the dual basis, the Fine-Scale Greens' function from \eqref{eq:g_prime_mu} may be computed componentwise. We start with the term $\mathcal{G}\boldsymbol{\mu}^T$ which is a row vector defined as
\begin{gather}
    \mathcal{G}\boldsymbol{\mu}^T = [\mathcal{G}\mu_0, \mathcal{G}\mu_1, \hdots \mathcal{G}\mu_{k}], \quad \text{with } \mathcal{G}\mu_i := \int_{\Omega} g(x, s) \mu_i(s) \: \mathrm{d}{s},
\end{gather}
where the components of the vector correspond to the solution of the underlying PDE (the Poisson equation in this case) with $\boldsymbol{\mu}$ as right hand side term. Similarly, the matrix $\boldsymbol{\mu} \mathcal{G} \boldsymbol{\mu}^T$ is computed as follows 
\begin{gather}
    \boldsymbol{\mu} \mathcal{G}\boldsymbol{\mu}^T = \left[\begin{array}{ccc}
        \prescript{}{V^*}{\left<\mu_0, \mathcal{G}\mu_0 \right>}_{V} & \hdots & \prescript{}{V^*}{\left<\mu_0, \mathcal{G}\mu_k \right>}_{V} \\
        \vdots & \ddots & \vdots \\
        \prescript{}{V^*}{\left<\mu_k, \mathcal{G}\mu_0 \right>}_{V} & \hdots & \prescript{}{V^*}{\left<\mu_k, \mathcal{G}\mu_k \right>}_{V}
    \end{array} \right],
    \label{eq:muGmu}
\end{gather}
where the columns can be recognised as the projection of $\mathcal{G} \boldsymbol{\mu}^T$ onto $\bar{V}$. As such, we have to respect that the entries of the matrix are duality pairings, meaning that for the $H_0^1$ projection, we have:
\begin{gather}
    \boldsymbol{\mu} \mathcal{G}\boldsymbol{\mu}^T = \left[\begin{array}{ccc}
        \left(\mu_0, \mathcal{G}\mu_0 \right)_{H_0^1(\Omega)} & \hdots & \left(\mu_0, \mathcal{G}\mu_k \right)_{H_0^1(\Omega)} \\
        \vdots & \ddots & \vdots \\
        \left(\mu_k, \mathcal{G}\mu_0 \right)_{H_0^1(\Omega)} & \hdots & \left(\mu_k, \mathcal{G}\mu_k \right)_{H_0^1(\Omega)}
    \end{array} \right], \quad \text{with }
    \left(\mu_i, \mathcal{G}\mu_j \right)_{H_0^1(\Omega)} := \int_{\Omega} \parddx{\mu_i}{x} \parddx{\mathcal{G}\mu_j}{x} \: \mathrm{d} x,
    \label{eq:muGmu_H01}
\end{gather}
and similarly for the $L^2$ projection we have:
\begin{gather}
    \boldsymbol{\mu} \mathcal{G}\boldsymbol{\mu}^T = \left[\begin{array}{ccc}
        \left(\mu_0, \mathcal{G}\mu_0 \right)_{L^2(\Omega)} & \hdots & \left(\mu_0, \mathcal{G}\mu_k \right)_{L^2(\Omega)} \\
        \vdots & \ddots & \vdots \\
        \left(\mu_k, \mathcal{G}\mu_0 \right)_{L^2(\Omega)} & \hdots & \left(\mu_k, \mathcal{G}\mu_k \right)_{L^2(\Omega)}
    \end{array} \right], \quad \text{with }
    \left(\mu_i, \mathcal{G}\mu_j \right)_{L^2(\Omega)} := \int_{\Omega} \mu_i \mathcal{G}\mu_j \: \mathrm{d} x.
    \label{eq:muGmu_L2}
\end{gather}

\begin{figure}[H]
\begin{subfigure}{0.49\linewidth}
    \centering
    \includegraphics[width = \linewidth]{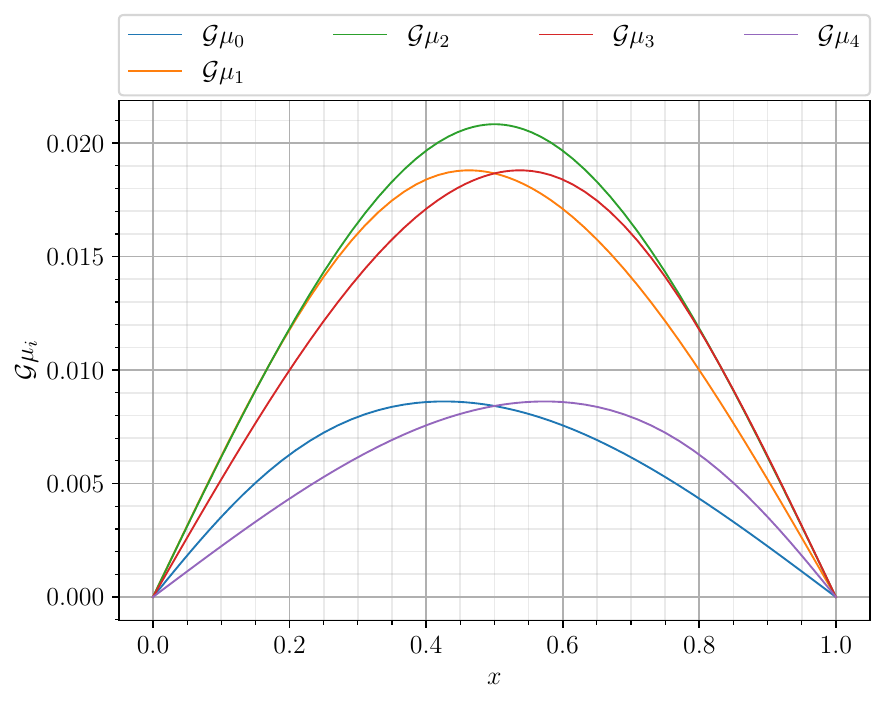}
    \caption{$\mathcal{G} \mu_i$ with $H_0^1$ dual basis functions from \Cref{fig:H01_dual_p3} as source term}
\end{subfigure}
\begin{subfigure}{0.49\linewidth}
    \centering
    \includegraphics[width = \linewidth]{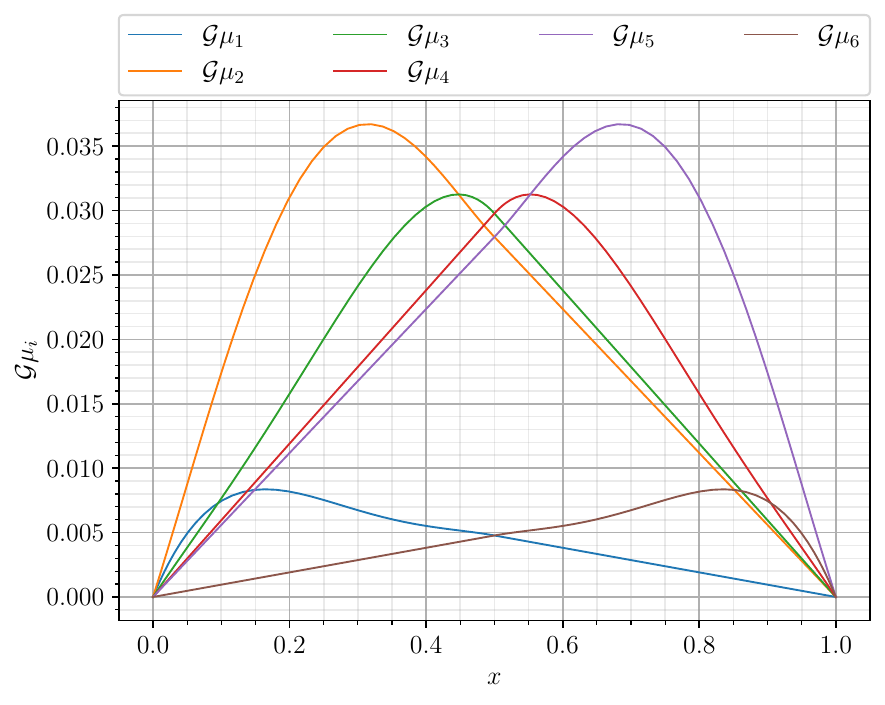}
    \caption{$\mathcal{G} \mu_i$ with $L^2$ dual basis functions from \Cref{fig:L2_dual_p2} as source term}
\end{subfigure}
\caption{$\mathcal{G} \mu_i$ computed with $H_0^1$ and $L^2$ dual basis functions for the Poisson equation}
\label{fig:mu_G_mu_L2_H01}
\end{figure}

The components of $\mathcal{G}\boldsymbol{\mu}^T$ are plotted in \Cref{fig:mu_G_mu_L2_H01} for both the $H_0^1$ and $L^2$ projections. $\boldsymbol{\mu} \mathcal{G}\boldsymbol{\mu}^T$ is then found by using \eqref{eq:muGmu} with the computed $\mathcal{G}\boldsymbol{\mu}^T$.

\begin{figure}[H]
\begin{subfigure}{0.49\linewidth}
    \centering
    \includegraphics[width = \linewidth]{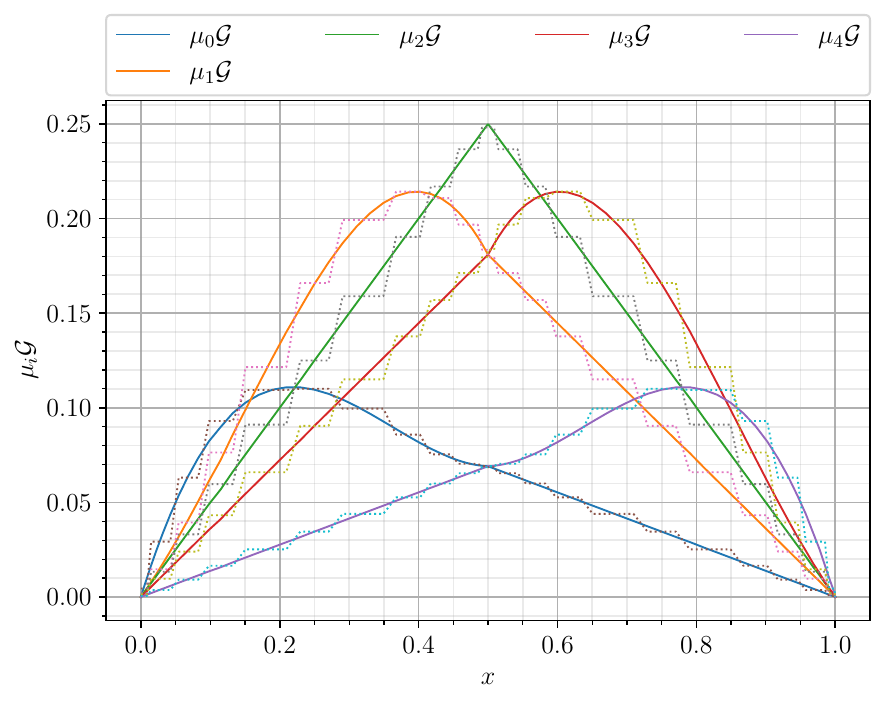}
    \caption{$\boldsymbol{\mu} \mathcal{G}$ computed using the $H_0^1$ dual basis functions from \Cref{fig:H01_dual_p3}}
    \label{fig:muG_3_H01}
\end{subfigure}
\begin{subfigure}{0.49\linewidth}
    \centering
    \includegraphics[width = \linewidth]{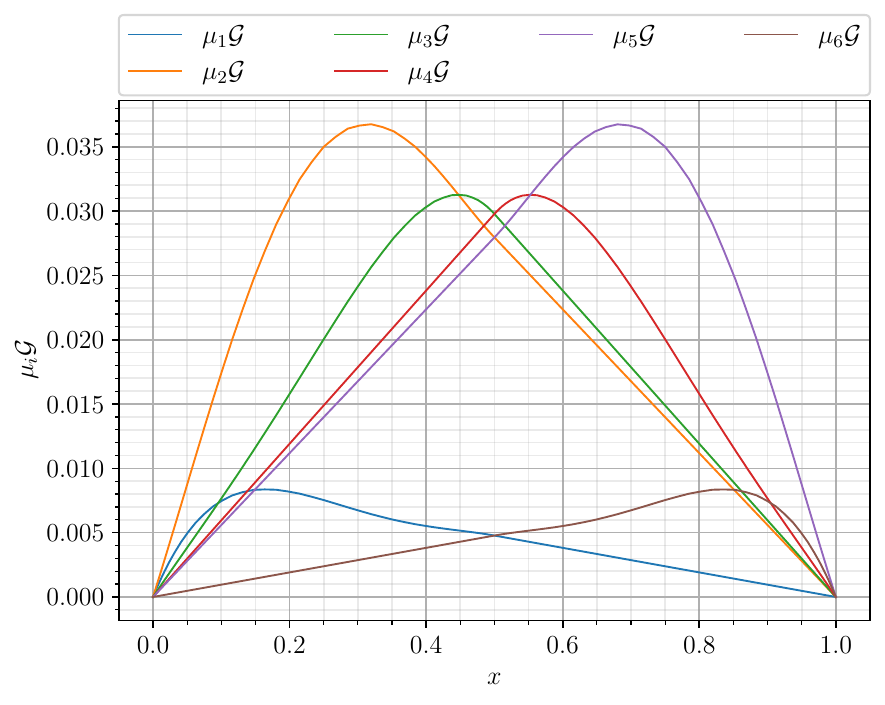}
    \caption{$\boldsymbol{\mu} \mathcal{G}$ computed using the $L^2$ dual basis functions from \Cref{fig:L2_dual_p2}}
    \label{fig:muG_3_L2}
\end{subfigure}
\caption{$\mu_i \mathcal{G}$ computed with $H_0^1$ and $L^2$ dual basis functions for the Poisson equation with the solid and dotted lines in \Cref{fig:muG_3_H01} corresponding to the discontinuity splitting integral approach and the standard quadrature rule, respectively}
\label{fig:muG}
\end{figure}

Lastly, we consider the term $\boldsymbol{\mu} \mathcal{G}$ which is a functional expressed as a column vector defined as
\begin{gather}
    \boldsymbol{\mu} \mathcal{G}(\nu) = \left[\begin{array}{c}
        \prescript{}{V^*}{\left<\mu_0, \mathcal{G}\nu \right>}_{V} \\
        \vdots \\
        \prescript{}{V^*}{\left<\mu_k, \mathcal{G}\nu \right>}_{V}
    \end{array} \right], \quad \forall \nu \in V^*.
\end{gather}
The act of applying $\boldsymbol{\mu} \mathcal{G}$ to an arbitrary function $\nu \in V^*$ is equivalent to computing the projection of the solution to the underlying PDE with $\nu$ as its source term. Computing the components of $\boldsymbol{\mu} \mathcal{G}$ is rather complicated, specifically for the $H_0^1$ projection. In particular, the issue arises from the duality pairing which requires the derivative of the piece-wise smooth Greens' function
\begin{gather}
    \left(\mu_i, \mathcal{G}\nu \right)_{H_0^1(\Omega)} = \left[\int_{\Omega} \parddx{\mu_i}{x} \parddx{\mathcal{G}}{x} \: \mathrm{d}x \right] \: \nu.
\end{gather}
Attempting to compute these integrals using standard quadrature rules will yield very poor results stemming from the fact that the discontinuity in $\parddx{\mathcal{G}}{x}$ will not be adequately captured by the quadrature. To resolve this issue, we split up the integral over the domain based on the known location of the discontinuity of $\parddx{\mathcal{G}}{x}$ (along $x = s$) and use the standard Gauss Lobatto quadrature to integrate the smooth functions on either side of the discontinuity. In contrast, for the $L^2$ case, $\boldsymbol{\mu} \mathcal{G}$ may be computed directly with standard quadrature as it only involves the $L^2$ pairing between $\boldsymbol{\mu}$ and $\mathcal{G}$. Plots for the components of $\boldsymbol{\mu} \mathcal{G}$ computed using the dual basis functions from \Cref{fig:H01_dual_p3} and \Cref{fig:L2_dual_p2} are shown in \Cref{fig:muG}. Moreover, \Cref{fig:muG_3_H01} shows $\boldsymbol{\mu} \mathcal{G}$ computed using two separate methods, where the solid lines correspond to the aforementioned discontinuity splitting integral approach and the dotted lines corresponding to the standard quadrature rule 
In the specific case of the Poisson equation and the $H_0^1$ projection, the entries of $\boldsymbol{\mu} \mathcal{G}$ simply yield the $\boldsymbol{\mu}$'s themselves as is apparent in \Cref{fig:muG_3_H01}. This is attributed to the fact that the $H_0^1$ projection encodes the Poisson problem as the projector is the differential operator. Consider $u_{exact}$ to be the exact (strong) solution to \eqref{eq:poisson_eq} satisfying \eqref{eq:poisson_bc}, from \eqref{eq:u_H01_proj_dual} we note that the $H_0^1$ projection ($\bar{u}$) of $u_{exact}$ is given by 
\begin{gather}
    \mathcal{N}^0(\bar{u}) = \left(\boldsymbol{\mu}^{H^1_0}, u_{exact} \right)_{H^1_0(\Omega)} \quad \xrightarrow[]{} \quad \bar{u} = \boldsymbol{\psi}^{(0)} \mathcal{N}^0(\bar{u}).
\end{gather}
Applying the pairing and employing integration by parts yields
\begin{gather}
    \mathcal{N}^0(\bar{u}) = \int_{\Omega} \parddx{\boldsymbol{\mu}^{H^1_0}}{x} \parddx{u_{exact}}{x} \: \mathrm{d}x = \cancel{\boldsymbol{\mu}^{H^1_0} \parddx{u_{exact}}{x} \Bigg|_{\partial \Omega}} - \int_{\Omega} \boldsymbol{\mu}^{H^1_0} \parddxdx{u_{exact}}{x} \: \mathrm{d}x,
\end{gather}
where the boundary term emerging from the integration by parts vanishes due to $\boldsymbol{\mu}^{H^1_0}$ being zero at the boundaries. Moreover, since $u_{exact}$ satisfies \eqref{eq:poisson_eq}, we have
\begin{equation}
    \mathcal{N}^0(\bar{u}) =  \int_{\Omega} \boldsymbol{\mu}^{H^1_0} f \: \mathrm{d}x,
\end{equation}
by which is it evident that the Poisson problem is embedded in the $H_0^1$ projection. It further implies that the exact $H_0^1$ projection of the Poisson equation can be computed without ever requiring the exact solution.

\begin{figure}[H]
\begin{subfigure}{0.49\linewidth}
    \centering
    \includegraphics[width = \linewidth]{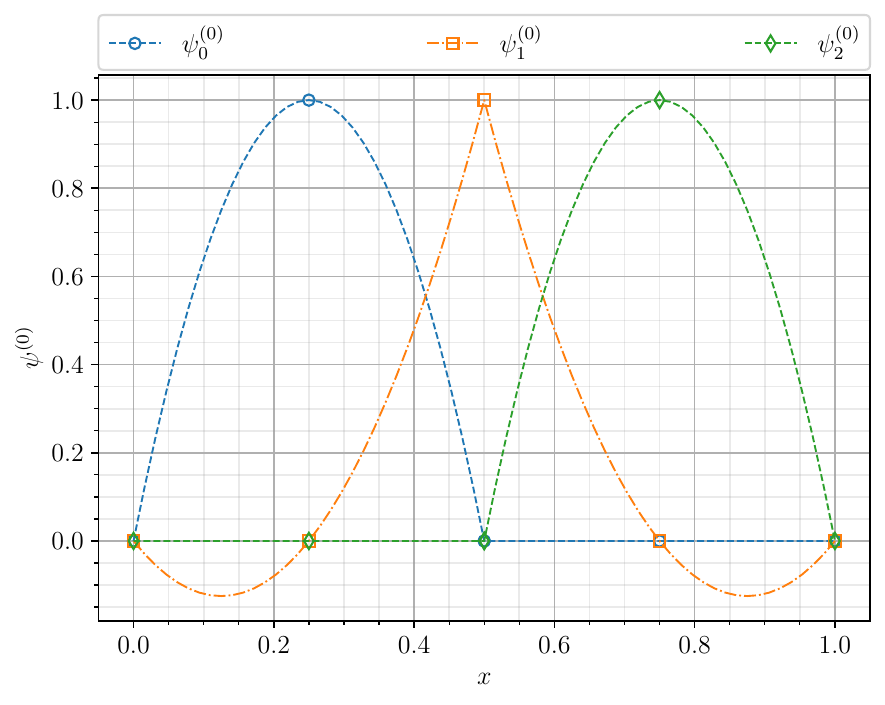}
    \caption{Nodal basis}
\end{subfigure}
\begin{subfigure}{0.49\linewidth}
    \centering
    \includegraphics[width = \linewidth]{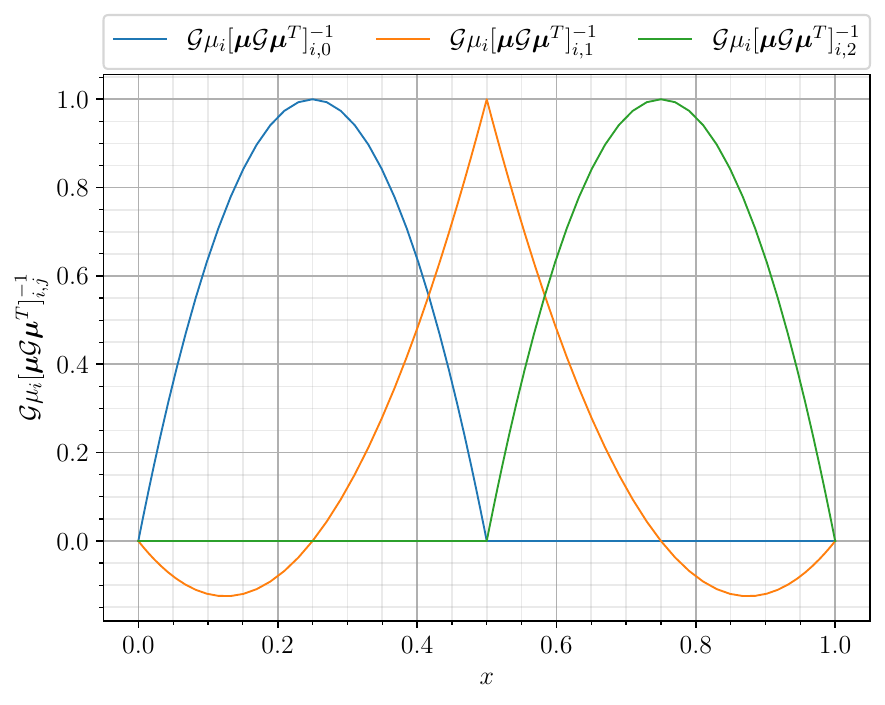}
    \caption{$\mathcal{G}\boldsymbol{\mu}^T [\boldsymbol{\mu} \mathcal{G}\boldsymbol{\mu}^T]^{-1}$ with $\boldsymbol{\mu} = \boldsymbol{\mu}^{H^1_0}$}
\end{subfigure}
\caption{Mesh with polynomial degree $p = 2$ and 2 elements}
\label{fig:Gmu_muGmu_H01_p2}
\end{figure}

\begin{figure}[H]
\begin{subfigure}{0.49\linewidth}
    \centering
    \includegraphics[width = \linewidth]{Images/edge_basis_N_2_p_3.pdf}
    \caption{Edge basis}
\end{subfigure}
\begin{subfigure}{0.49\linewidth}
    \centering
    \includegraphics[width = \linewidth]{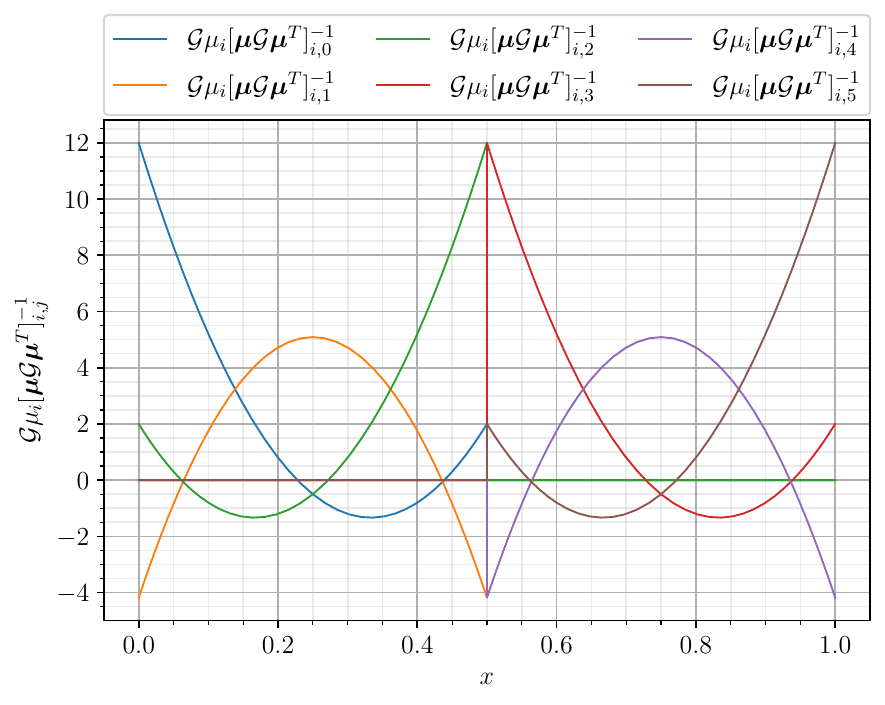}
    \caption{$\mathcal{G}\boldsymbol{\mu}^T [\boldsymbol{\mu} \mathcal{G}\boldsymbol{\mu}^T]^{-1}$ with $\boldsymbol{\mu} = \boldsymbol{\mu}^{L^2}$}
\end{subfigure}
\caption{Mesh with polynomial degree $p = 3$ and 2 elements}
\label{fig:Gmu_muGmu_L2_p2}
\end{figure}

With all these components in place, the Fine-Scale Greens' function may be assembled. Before doing so, however, it is worthwhile reviewing the sequence of mappings associated with the Fine-Scale Greens' function. The full expression for computing $u'$ reads
\begin{equation}
    u' = \mathcal{G}\mathscr{R} \bar{u} - \mathcal{G}\boldsymbol{\mu}^T [\boldsymbol{\mu} \mathcal{G}\boldsymbol{\mu}^T]^{-1} \boldsymbol{\mu} \mathcal{G}\mathscr{R} \bar{u}, \label{eq:u_prime_g_prime}
\end{equation}
where we see the first term ($\mathcal{G}\mathscr{R} \bar{u}$) is the (global) Greens' function applied to the residual of the resolved scales. The same term appears in the second term of \eqref{eq:u_prime_g_prime}, but there it is accompanied by $\boldsymbol{\mu}$, which indicates it is being projected onto the finite-dimensional resolved space. The remaining component $\mathcal{G}\boldsymbol{\mu}^T [\boldsymbol{\mu} \mathcal{G}\boldsymbol{\mu}^T]^{-1}$ is rather interesting as it reproduces the resolved scale basis as shown in \Cref{fig:Gmu_muGmu_H01_p2} and \Cref{fig:Gmu_muGmu_L2_p2}.

The term $\mathcal{G}\boldsymbol{\mu}^T [\boldsymbol{\mu} \mathcal{G}\boldsymbol{\mu}^T]^{-1}$ can thus be seen as a reconstruction/interpolation operator mapping discrete quantities from the resolved finite-dimensional space to the continuous space. Therefore, the term $\mathcal{G}\boldsymbol{\mu}^T [\boldsymbol{\mu} \mathcal{G}\boldsymbol{\mu}^T]^{-1} \boldsymbol{\mu} \mathcal{G}\mathscr{R} \bar{u}$ simply equates to the component of $\mathcal{G}\mathscr{R} \bar{u}$ that lives in the resolved space. Note, it is not a projection of the residual, but rather the projection of $\mathcal{G}$ applied to the residual.

\begin{figure}[H]
\begin{subfigure}{0.49\linewidth}
    \centering
    \includegraphics[width = \linewidth]{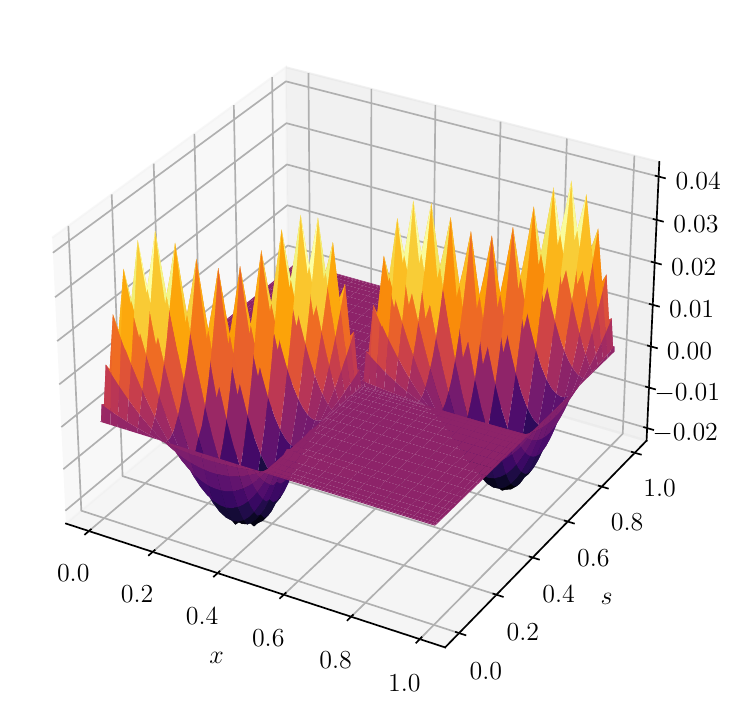}
    \caption{Fine-scale Greens' function for the $H_0^1$ projection}
    \label{fig:Gp_h01_p2}
\end{subfigure}
\begin{subfigure}{0.49\linewidth}
    \centering
    \includegraphics[width = \linewidth]{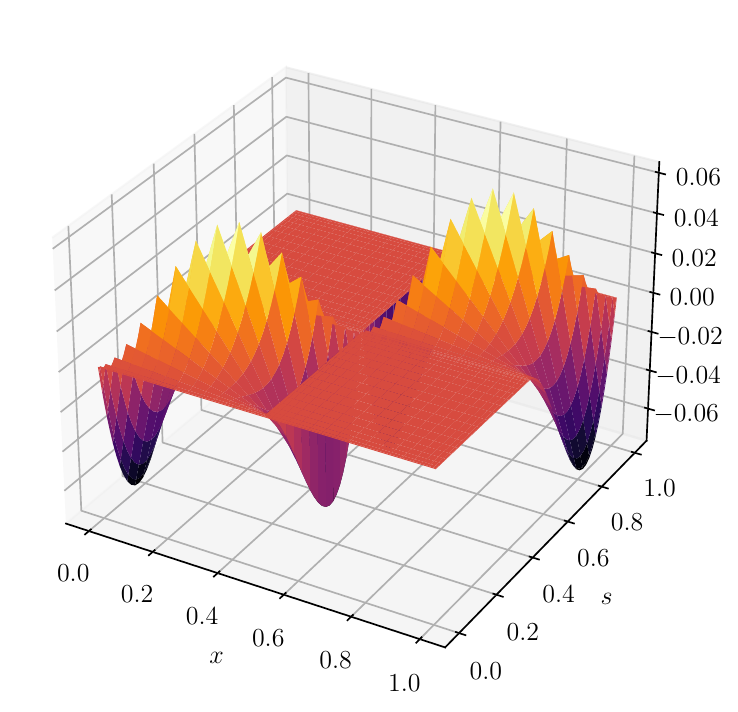}
    \caption{Fine-Scale Greens' function for the $L^2$ projection}
    \label{fig:Gp_L2_p2}
\end{subfigure}
\caption{Fine-Scale Greens' functions for 1D Poisson equation on a mesh with polynomial degree $p = 2$ and 2 elements}
\end{figure}

\begin{figure}[H]
\begin{subfigure}{0.49\linewidth}
    \centering
    \includegraphics[width = \linewidth]{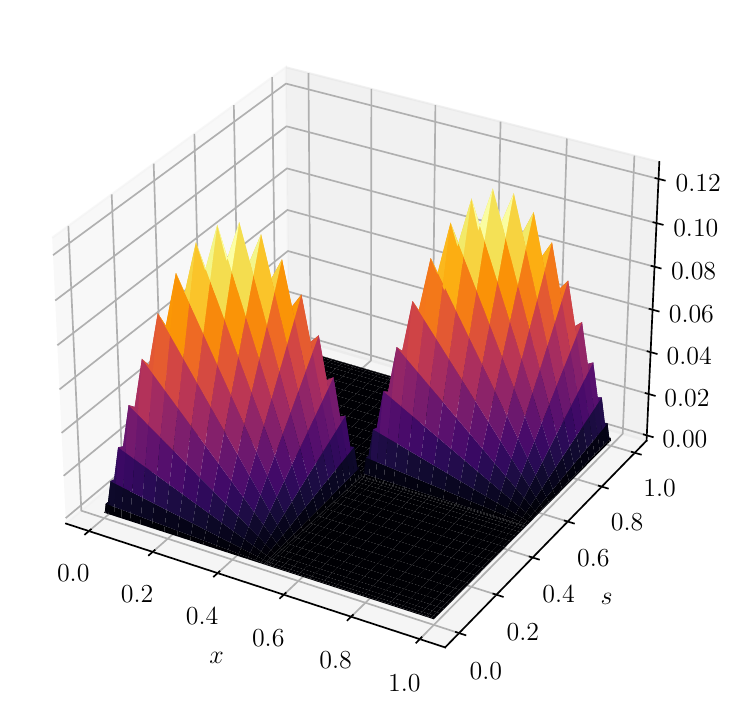}
    \caption{Element's Greens' functions}
    \label{fig:G_el_p1}
\end{subfigure}
\begin{subfigure}{0.49\linewidth}
    \centering
    \includegraphics[width = \linewidth]{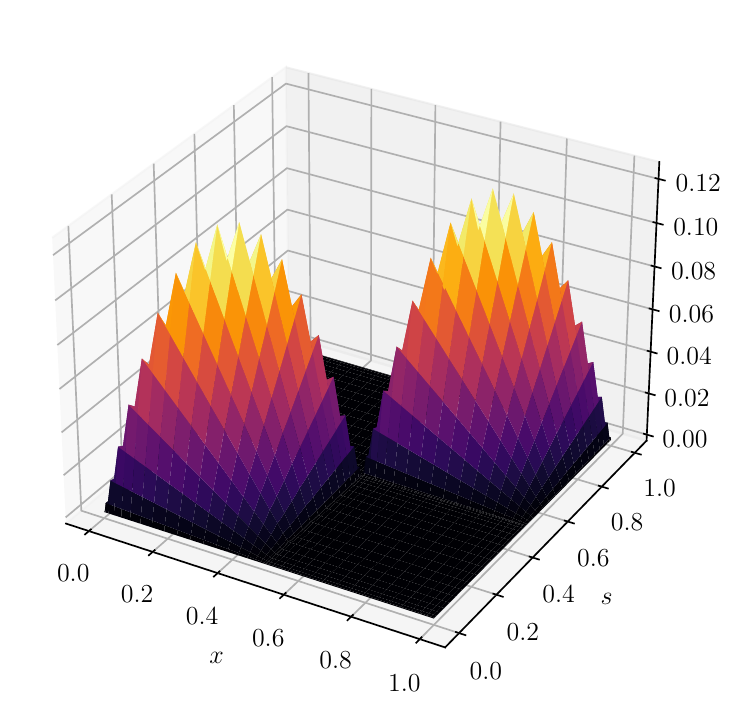}
    \caption{Fine-Scale Greens' function for the $H_0^1$ projection}
    \label{fig:Gp_p1}
\end{subfigure}
\caption{Greens' functions for 1D Poisson equation on a mesh with polynomial degree $p = 1$ and 2 elements}
\label{fig:Gp_Gel_p1}
\end{figure}

\begin{figure}[H]
\begin{subfigure}{0.49\linewidth}
    \centering
    \includegraphics[width = \linewidth]{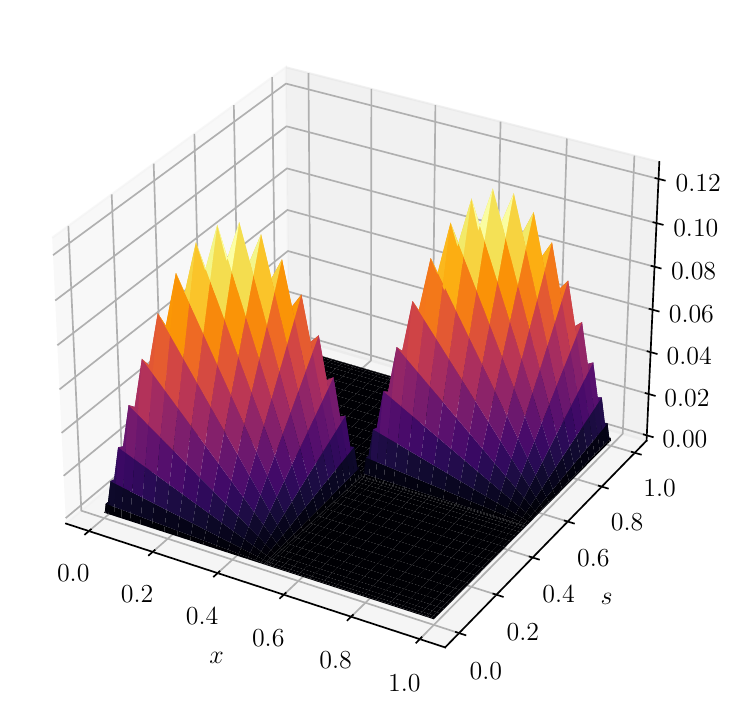}
    \caption{Element's Greens' functions}
    \label{fig:G_el_p2}
\end{subfigure}
\begin{subfigure}{0.49\linewidth}
    \centering
    \includegraphics[width = \linewidth]{Images/Gp_H01_N_2_p_2_c_0.000_nu_1.000.pdf}
    \caption{Fine-Scale Greens' function for the $H_0^1$ projection}
    \label{fig:Gp_p2}
\end{subfigure}
\caption{Greens' functions for 1D Poisson equation on a mesh with polynomial degree $p = 2$ and 2 elements}
\label{fig:Gp_Gel_p2}
\end{figure}

\Cref{fig:Gp_h01_p2} and \Cref{fig:Gp_L2_p2} show the resulting Fine-Scale Greens' functions for the $H_0^1$ and $L^2$ projections, respectively. Furthermore, \Cref{fig:Gp_Gel_p1} and \Cref{fig:Gp_Gel_p2} show the $H_0^1$ Fine-Scale Greens' function alongside the localised Element's Greens' function for $p = 1$ and $p = 2$. Here we find that the Fine-Scale Greens' function for the $p = 1$ case exactly corresponds to the local Element's Greens' function. This, however, is no longer true for $p > 1$.

\subsection{Reconstruction of fine-scale terms of a projection for 1D Poisson}

\begin{figure}[H]
\centering
\begin{subfigure}{0.49\linewidth}
    \centering
    \includegraphics[width = \linewidth]{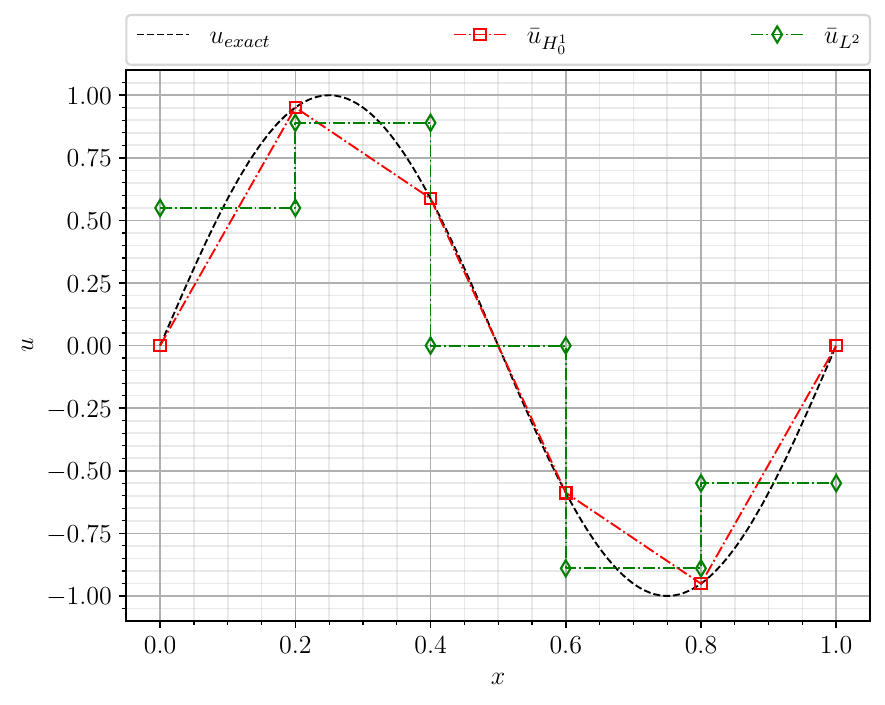}
    \caption{Exact solution to sample Poisson equation and its $H_0^1$ and $L^2$ projections}
    \label{fig:u_ex_proj_p1_N5}
\end{subfigure} \\
\begin{subfigure}{0.49\linewidth}
    \centering
    \includegraphics[width = \linewidth]{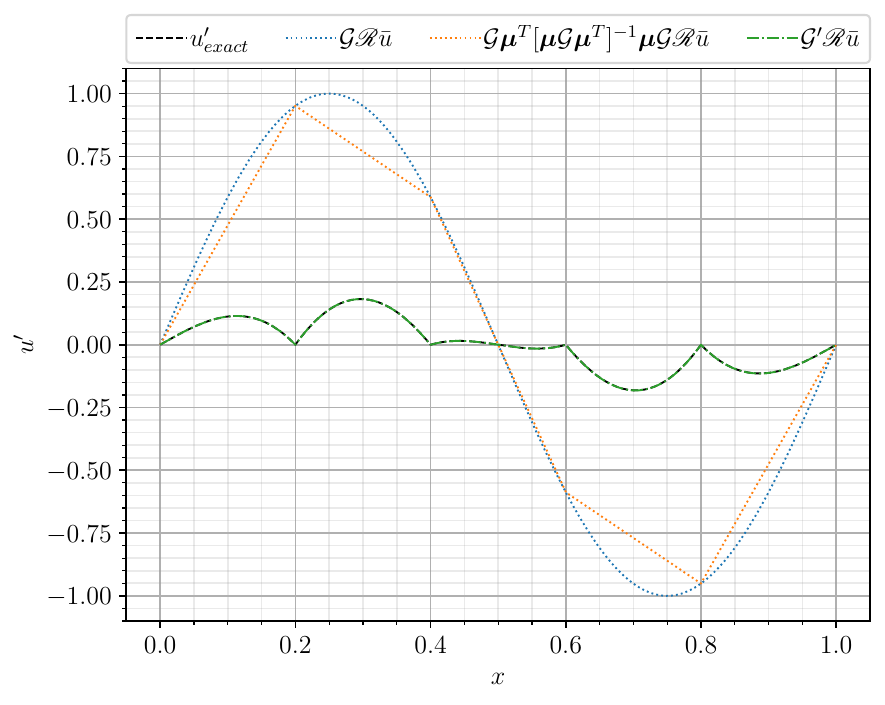}
    \caption{Fine scales for $H_0^1$ projection computed using the Fine-Scale Greens' functions}
    \label{fig:u_p_H01_p1_N5}
\end{subfigure}
\begin{subfigure}{0.49\linewidth}
    \centering
    \includegraphics[width = \linewidth]{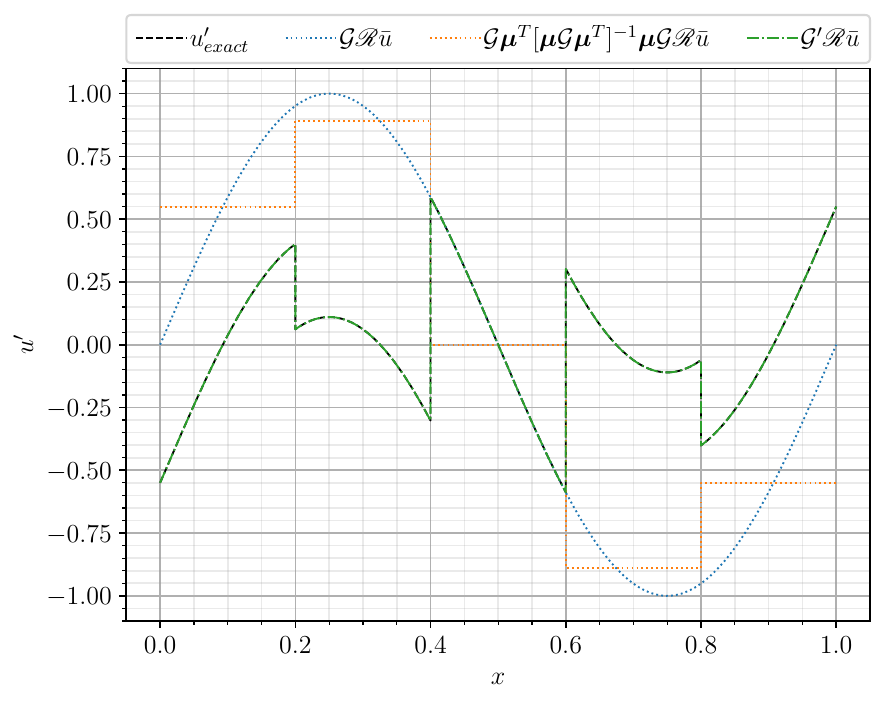}
    \caption{Fine scales for $L^2$ projection computed using the Fine-Scale Greens' functions}
     \label{fig:u_p_L2_p1_N5}
\end{subfigure}
\caption{Mesh with polynomial degree $p = 1$ and 5 elements}
\end{figure}

\begin{figure}[H]
\centering
\begin{subfigure}{0.49\linewidth}
    \centering
    \includegraphics[width = \linewidth]{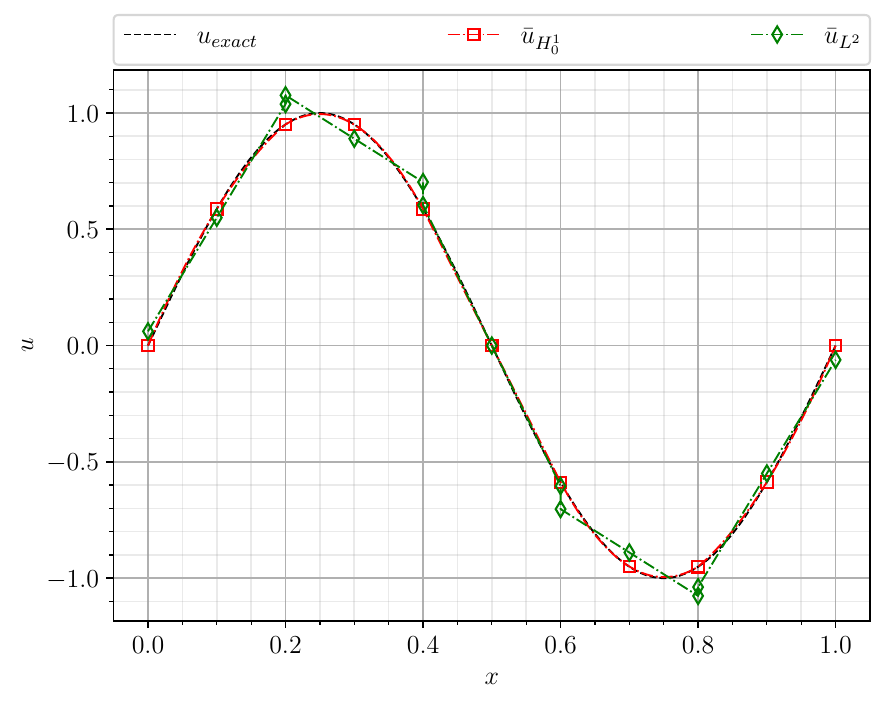}
    \caption{Exact solution to sample Poisson equation and its $H_0^1$ and $L^2$ projections}
    \label{fig:u_ex_proj_p2_N5}
\end{subfigure} \\
\begin{subfigure}{0.49\linewidth}
    \centering
    \includegraphics[width = \linewidth]{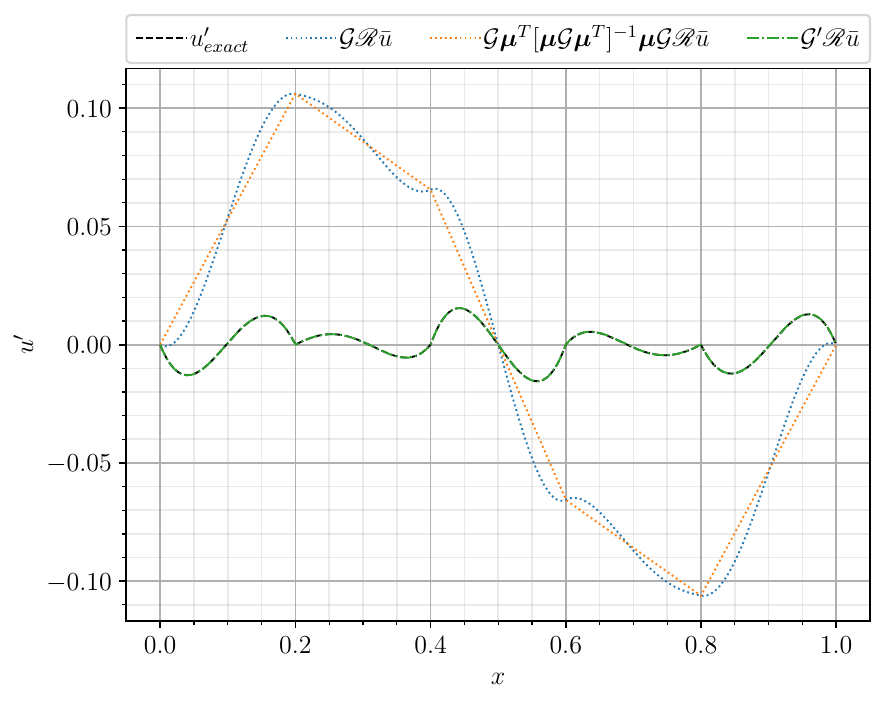}
    \caption{Fine scales for $H_0^1$ projection computed using the Fine-Scale Greens' functions}
     \label{fig:u_p_H01_p2_N5}
\end{subfigure}
\begin{subfigure}{0.49\linewidth}
    \centering
    \includegraphics[width = \linewidth]{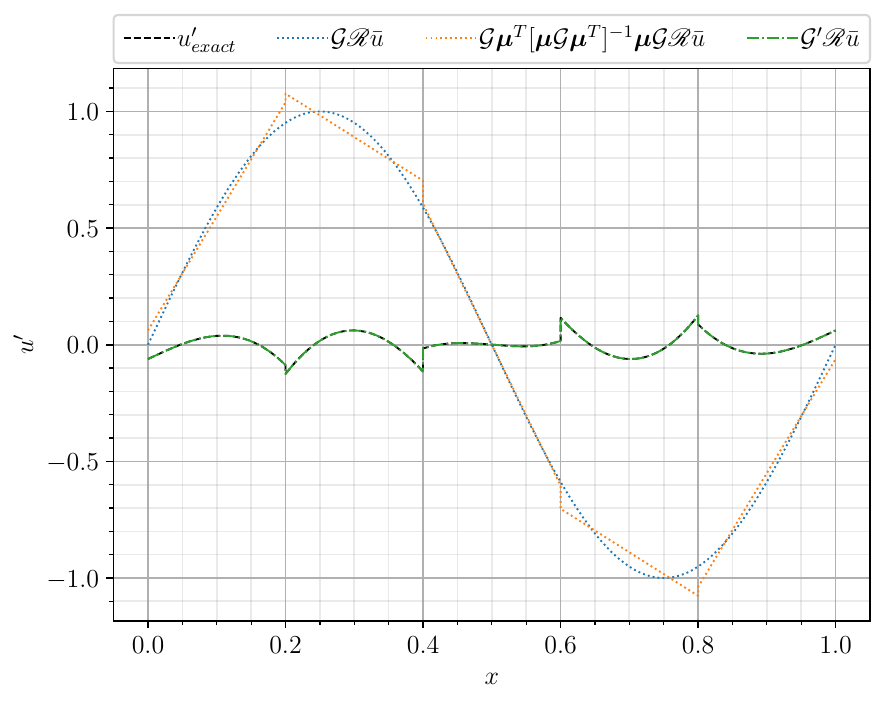}
    \caption{Fine scales for $L^2$ projection computed using the Fine-Scale Greens' functions}
     \label{fig:u_p_L2_p2_N5}
\end{subfigure}
\caption{Mesh with polynomial degree $p = 2$ and 5 elements}
\end{figure}

In order to carry out numerical tests with the constructed Fine-Scales Greens' function, we consider the sample Poisson problem from 
\eqref{eq:poisson_eq} with
\begin{equation}
    f = 4 \pi^2 \sin{(2 \pi x)},
\end{equation}
for which the exact solution is given by
\begin{equation}
    u_{exact} = \sin{(2 \pi x)}.
\end{equation}
We then compute the exact $H_0^1$ and $L^2$ projections of the exact solution using the corresponding dual basis functions as previously described. The plots of the exact solution along with its $H_0^1$ and $L^2$ projections onto the finite-dimensional mesh are shown in \Cref{fig:u_ex_proj_p1_N5} and \Cref{fig:u_ex_proj_p2_N5} for two polynomial approximations. Once again, we choose a $L^2$ projection that maps onto the space of edge polynomials following the line the reasoning in \Cref{rem:L2_space}. As a consequence, the degree of polynomials spanning the $L^2$ space is $(p - 1)$ (see \Cref{rem:edge_p}).

Applying the respective Fine-Scale Greens' functions to the residual of the respective projections yields the exact unresolved scales of the projections as shown in \Cref{fig:u_p_H01_p1_N5} and \Cref{fig:u_p_H01_p2_N5} for the $H_0^1$ projection and \Cref{fig:u_p_L2_p1_N5} and \Cref{fig:u_p_L2_p2_N5} for the $L^2$ projection. These figures also include the plots of the individual components $\mathcal{G} \mathscr{R} \bar{u}$ and $\mathcal{G}\boldsymbol{\mu}^T \left[\boldsymbol{\mu} \mathcal{G} \boldsymbol{\mu}^T \right]^{-1} \boldsymbol{\mu} \mathcal{G} \mathscr{R}\bar{u}$. These individual components highlight the main characteristic of the Fine-Scale Greens' function, namely the subtraction of the component of $\mathcal{G} \mathscr{R} \bar{u}$ that lives in the resolved space. 

\subsection{Computing the fine-scale terms of a projection for 1D advection-diffusion}

In a manner similar to that considered for the 1D Poisson equation, we can apply this approach to construct the Fine-Scale Greens' function for the 1D advection-diffusion equation. The governing equation for the 1D advection-diffusion problem reads as follows
\begin{gather}
    c \parddx{u}{x} - \nu \parddxdx{u}{x} = f, \quad x \in \Omega = [0, 1] \label{eq:adv_diff}\\
    u(x) = 0, \quad x \in \partial \Omega,
\end{gather}
with the constants, $c$ and $\nu$ being the advection speed and diffusion coefficient, respectively.

In this simple 1D setting, the global Greens' function for the advection-diffusion equation in $\Omega$ can be found to be
\begin{gather}
    g(x, s) = \begin{cases}
        \frac{1 - e^{-2 \alpha \left(1 - \frac{x}{h} \right)}}{c \left(1 - e^{-2 \alpha} \right)} \left(1 - e^{-2 \alpha \frac{s}{h}} \right), \quad \: \: \quad &x \leq s \\
        \frac{e^{2 \alpha \frac{x}{h}} - 1}{c \left(1 - e^{-2 \alpha} \right)} \left(e^{-2 \alpha \frac{s}{h}} - e^{-2 \alpha} \right), \quad &x > s
    \end{cases} \label{eq:greens_ad} \\
    \alpha = \frac{c h}{2 \nu}, \label{def:alpha}
\end{gather}
where $h$ is the domain width and $\alpha$ is the mesh Peclet number. In theory, we can take this Greens' function and construct the Fine-Scale Greens' function for any projector using the aforementioned approach. This poses two notable complications, however, one relating to the integration of the Greens' function's sharp gradients at high Peclet numbers and the second being the loss of generality in multi-dimensional settings where the global Greens' function of the advection-diffusion problem is not readily available. We therefore adopt an alternate approach for the advection-diffusion problem which bypasses these shortcomings. 

This entails rewriting the advection-diffusion equation as a Poisson problem with a modified right hand side term as follows
\begin{gather}
    -\parddxdx{u}{x} = \frac{1}{\nu} f - \frac{c}{\nu} \parddx{u}{x},
\end{gather}
for which the resolved component of the variational form reads
\begin{gather}
    -\prescript{}{V^*}{\left\langle \parddxdx{\bar{u}}{x}, \Bar{v} \right\rangle_{V}} - \prescript{}{V^*}{\left\langle \parddxdx{u'}{x}, \Bar{v} \right\rangle_{V}} = \prescript{}{V^*}{\left\langle \left(\frac{1}{\nu} f - \frac{c}{\nu} \parddx{u}{x} \right), \Bar{v} \right\rangle_{V}} \quad \forall \Bar{v} \in \Bar{V}, \label{eq:AD_new}
\end{gather}
and the fine scales may be computed using the Fine-Scale Greens' function for Poisson equation as follows
\begin{gather}
    u' = \mathcal{G}' \mathscr{R} \bar{u} = \mathcal{G}' \left(\left(\frac{1}{\nu} f - \frac{c}{\nu} \parddx{u}{x}\right) + \parddxdx{\bar{u}}{x} \right).
\end{gather}
Note that in the diffusive limit with $c = 0$ this reverts back to the formulation for the Poisson equation. Another important thing to note is that the advection term $\parddx{u}{x}$ appearing in the residual is the gradient of the exact solution and not $\bar{u}$. This term thus needs special treatment which will be discussed shortly. We will first demonstrate how this formulation reconstructs the fine scales for the projections in \Cref{subsec:AD_reconst}, when we substitute the known exact derivative of the solution in the residual. Thereafter in \Cref{subsec:iterative}, we describe an iterative approach which overcomes the need to have the exact derivative of the solution. 

\subsection{Reconstruction of fine-scale terms of a projection for 1D advection-diffusion}
\label{subsec:AD_reconst}
\begin{figure}[H]
\centering
\begin{subfigure}{0.49\linewidth}
    \centering
    \includegraphics[width = \linewidth]{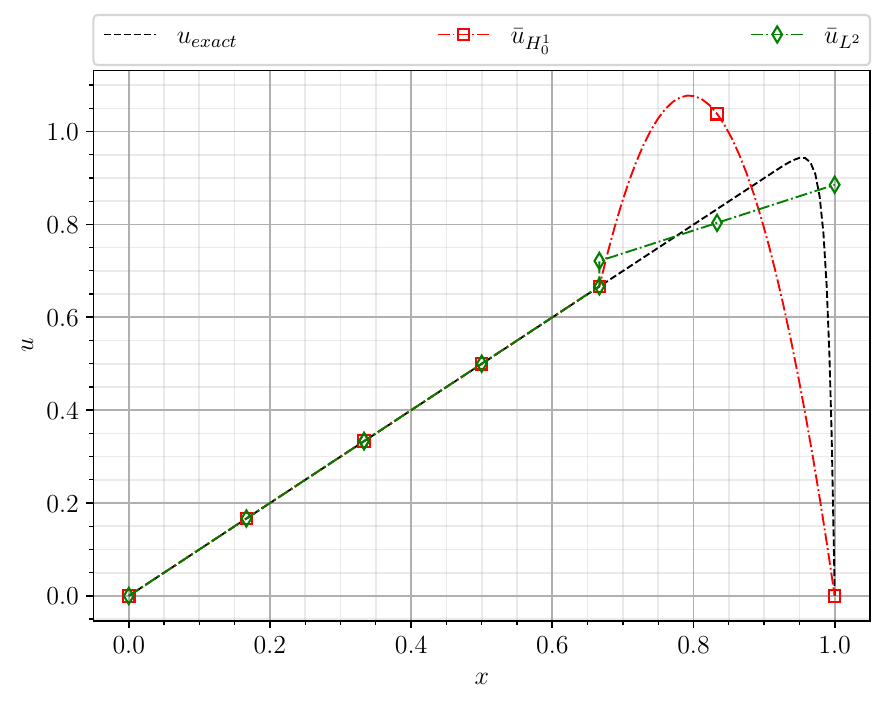}
    \caption{Exact solution to advection-diffusion equation and its $H_0^1$ and $L^2$ projections with $\alpha = 50$}
    \label{fig:u_ex_proj_ad1}
\end{subfigure} \\
\begin{subfigure}{0.49\linewidth}
    \centering
    \includegraphics[width = \linewidth]{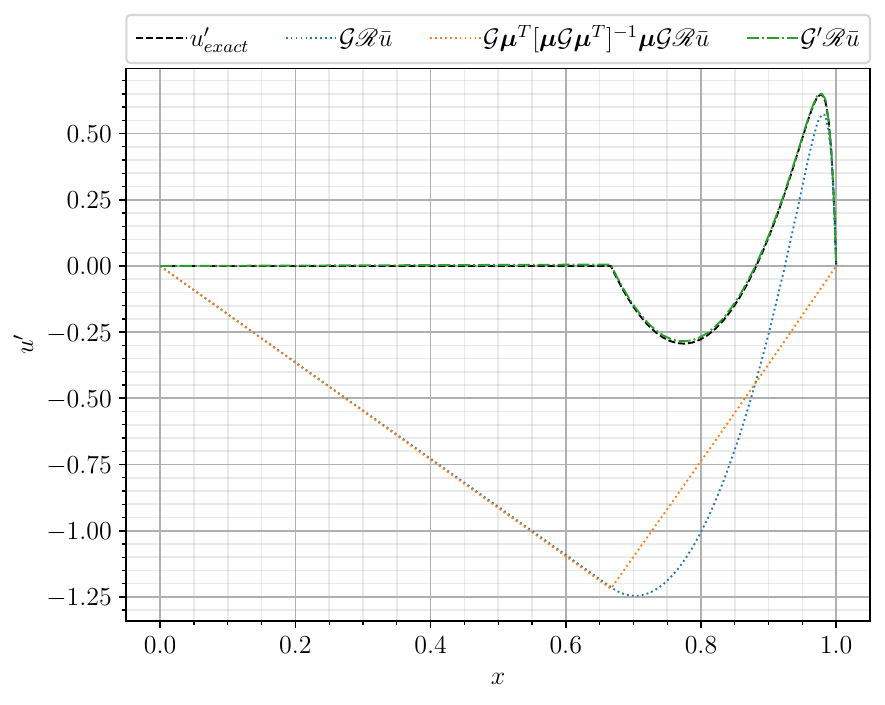}
    \caption{Fine-scales for $H_0^1$ projection computed using the Fine-Scale Greens' functions}
    \label{fig:u_p_H01_1}
\end{subfigure}
\begin{subfigure}{0.49\linewidth}
    \centering
    \includegraphics[width = \linewidth]{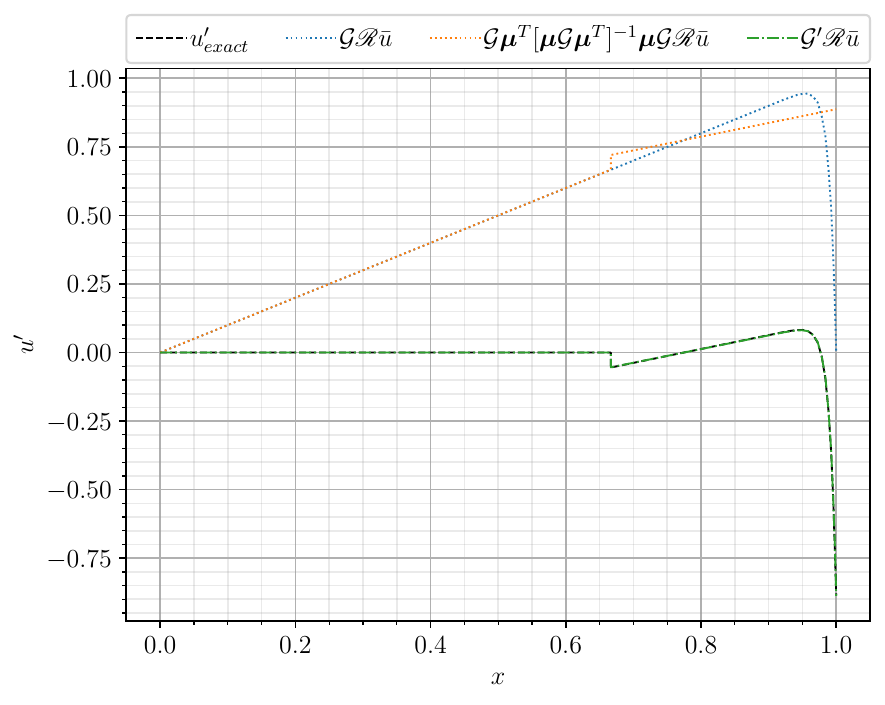}
    \caption{Fine-scales for $L^2$ projection computed using the Fine-Scale Greens' functions}
    \label{fig:u_p_L2_1}
\end{subfigure}
\caption{Mesh with polynomial degree $p = 2$ and 3 elements}
\end{figure}

\begin{figure}[H]
\centering
\begin{subfigure}{0.49\linewidth}
    \centering
    \includegraphics[width = \linewidth]{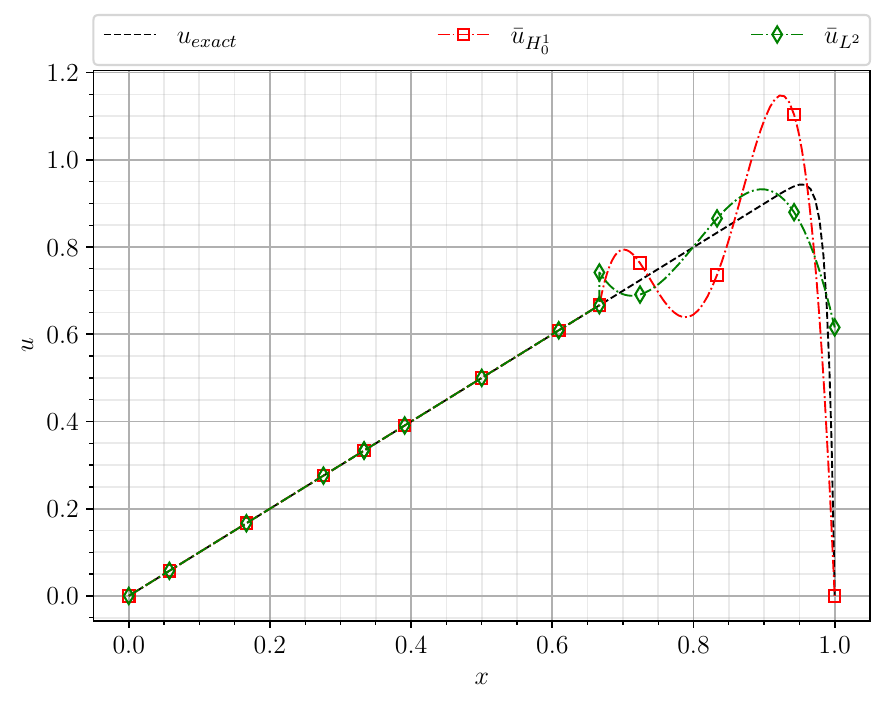}
    \caption{Exact solution to advection-diffusion equation and its $H_0^1$ and $L^2$ projections with $\alpha = 50$}
    \label{fig:u_ex_proj_ad2}
\end{subfigure} \\
\begin{subfigure}{0.49\linewidth}
    \centering
    \includegraphics[width = \linewidth]{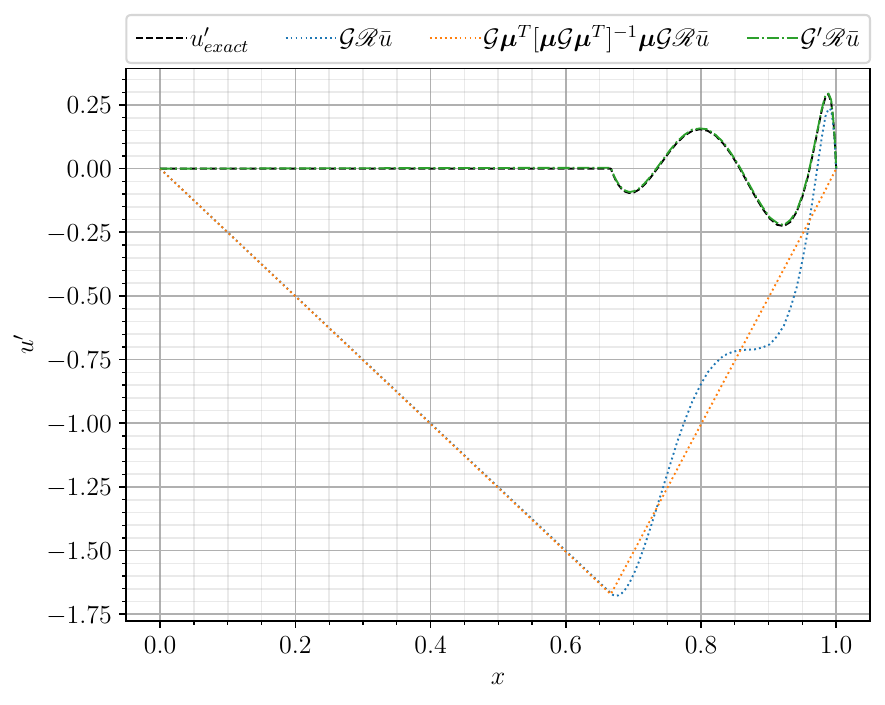}
    \caption{Fine-scales for $H_0^1$ projection computed using the Fine-Scale Greens' functions}
    \label{fig:u_p_H01_2}
\end{subfigure}
\begin{subfigure}{0.49\linewidth}
    \centering
    \includegraphics[width = \linewidth]{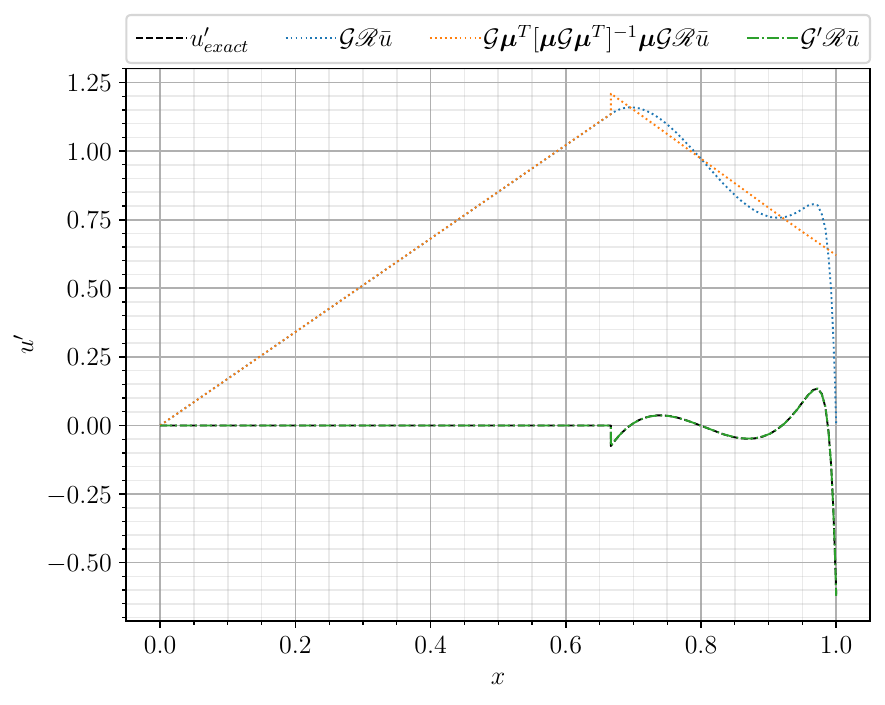}
    \caption{Fine-scales for $L^2$ projection computed using the Fine-Scale Greens' functions}
    \label{fig:u_p_L2_2}
\end{subfigure}
\caption{Mesh with polynomial degree $p = 4$ and 3 elements}
\end{figure}

We consider the advection-diffusion equation from \eqref{eq:adv_diff} with the following input parameters
\begin{equation}
    f = 1, \quad c = 1, \quad \nu = 0.01, \quad \alpha = 50 
\end{equation}
for which the exact solution and the exact solution gradient are given by
\begin{gather}
    u_{exact} = \frac{1}{c} \left(x - \frac{e^{\alpha(x - 1)} - e^{-\alpha}}{1 - e^{-\alpha}} \right) \\
    \parddx{u_{exact}}{x} = \frac{1}{c} \left(1 - \frac{\alpha e^{\alpha(x - 1)}}{1 - e^{-\alpha}} \right)
\end{gather}
The plot of the exact solution and its $H_0^1$ and $L^2$ projections onto different meshes are shown in \Cref{fig:u_ex_proj_ad1} and \Cref{fig:u_ex_proj_ad2}. The exact fine scales and the reconstruction thereof using the (Poisson) Fine-Scale Greens' function are shown in \Cref{fig:u_p_H01_1} and \Cref{fig:u_p_H01_2} for the $H_0^1$ projection and in \Cref{fig:u_p_L2_1} and \Cref{fig:u_p_L2_2} for the $L^2$ projection. As evident from these figures, the proposed formulation correctly reconstructs all the fine scales missing in the projection. However, as noted before, this approach required the gradient of the exact solution to be available which makes it unusable for any VMS formulation in practice. 

\subsection{Iterative VMS approach}
\label{subsec:iterative}
The need for the derivative of the exact solution in the above demonstration can be avoided using an iterative procedure. The goal of this iterative approach is to produce a VMS formulation for the advection-diffusion equation which yields the exact projection and the missing fine scales. Here, we focus on the $H_0^1$ projection, motivated by the fact that the exact solution lives in the $H_0^1$ space. We note, however, that this iterative concept is not limited to $H_0^1$ and can be extended to other projections.

To construct an iterative approach that solves for $H_0^1$ projection, we start with the resolved scale equation from \eqref{eq:AD_new} where we take $v = {\boldsymbol{\mu}}^{H_0^1}$ as the test functions and apply an $L^2$ pairing
\begin{gather}
    -\left(\boldsymbol{\mu}^{H_0^1}, \parddxdx{u}{x} \right)_{L^2(\Omega)} = \left(\boldsymbol{\mu}^{H_0^1}, \frac{1}{\nu} f \right)_{L^2(\Omega)} - \left(\boldsymbol{\mu}^{H_0^1}, \frac{c}{\nu} \parddx{u}{x} \right)_{L^2(\Omega)}.
\end{gather}
If we apply integration by parts to the left-hand side, we get
\begin{gather}
    \cancel{-\boldsymbol{\mu}^{H_0^1} \parddx{u}{x} \Bigg|_{\partial \Omega}} + \left(\parddx{\boldsymbol{\mu}^{H_0^1}}{x}, \parddx{u}{x} \right)_{L^2(\Omega)} = \left(\boldsymbol{\mu}^{H_0^1}, \frac{1}{\nu} f \right)_{L^2(\Omega)} - \cancel{\boldsymbol{\mu}^{H_0^1} u \Bigg|_{\partial \Omega}} + \left(\parddx{\boldsymbol{\mu}^{H_0^1}}{x}, \frac{c}{\nu} u \right)_{L^2(\Omega)},
\end{gather}
where the boundary terms cancel due to $\boldsymbol{\mu}^{H_0^1}$ being zero at the boundaries. Furthermore, the remaining term on the left can be identified to be a duality pairing in $H_0^1$, which exactly yields the expansion coefficients of the projection ($\mathcal{N}^0(\bar{u})$) as noted in \eqref{eq:u_H01_proj_dual}. We thus have the following equations for the projection/resolved scales and the fine scales
\begin{gather}
    \mathcal{N}^0(\bar{u}) = \left(\boldsymbol{\mu}^{H_0^1}, \frac{1}{\nu} f \right)_{L^2(\Omega)} + \left(\parddx{\boldsymbol{\mu}^{H_0^1}}{x}, \frac{c}{\nu}\left(\bar{u} + u'\right) \right)_{L^2(\Omega)} \quad \rightarrow \quad \quad \bar{u} = \boldsymbol{\psi}^{(0)} \mathcal{N}^0(\bar{u}) \\
    u' = \mathcal{G}' \mathscr{R} \bar{u} = \mathcal{G}' \left(\frac{1}{\nu} f - \frac{c}{\nu} \parddx{\bar{u}}{x} - \frac{c}{\nu} \parddx{u'}{x} + \parddxdx{\bar{u}}{x} \right), \label{eq:up_ad_iter}
\end{gather}
where use the fact that $u = \bar{u} + u'$. An important thing to note is that this formulation solves for the strong $u'$ which we choose to represent on a finer mesh\footnote{The finer mesh is solely used to represent/plot $u'$}. For conciseness, we will use the following shorthand notation to represent the above equations
\begin{align}
    \bar{u} &= \boldsymbol{\bar{\Phi}}(f, \bar{u}, u') \\
    u' &= \boldsymbol{\Phi'}(f, \bar{u}, u').
\end{align}
In order to solve this coupled set of equations, we employ an iterative scheme, where we initialise $\bar{u}_0 = 0$ and $u'_0 = 0$, and use the following update equations to compute $\bar{u}$ and $u'$ at the $i + 1$ iteration
\begin{align}
    \bar{u}_{i + 1} &= \bar{u}_{i} + w(\boldsymbol{\bar{\Phi}}(f, \bar{u}_i, u'_i) - \bar{u}_{i}) \\
    u'_{i + 1} &= u'_{i} + w (\boldsymbol{\Phi'}(f, \bar{u}_i, u'_i) - u'_{i}),
\end{align}
with $w$ as the under-relaxation factor which we set equal to $\frac{1}{2 \alpha}$, with $\alpha$ given by \eqref{def:alpha}.
We successively iterate the above equations until the residual drops below a user-specified tolerance $\varepsilon$, $\left\lVert\bar{u}_{i + 1} - \bar{u}_{i} \right\rVert_{L^2} < \varepsilon$. For all the results shown in this subsection, the tolerance was set to $\varepsilon = 10^{-8}$. 

\begin{figure}[H]
\begin{subfigure}{0.49\linewidth}
    \centering
    \includegraphics[width = \linewidth]{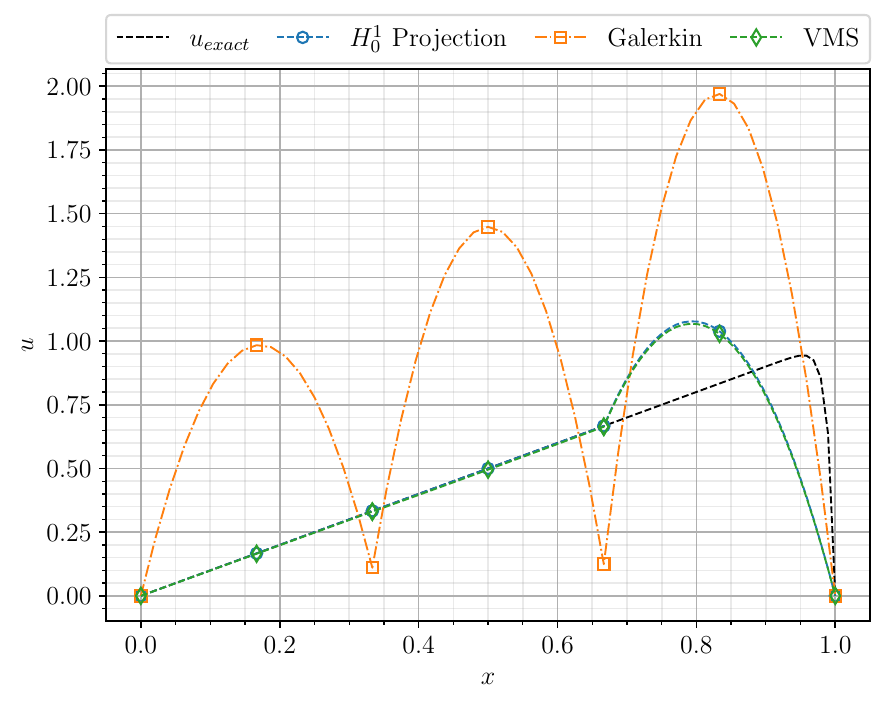}
    \caption{Polynomial degree $p = 2$ and 3 elements}
\end{subfigure}
\begin{subfigure}{0.49\linewidth}
    \centering
    \includegraphics[width = \linewidth]{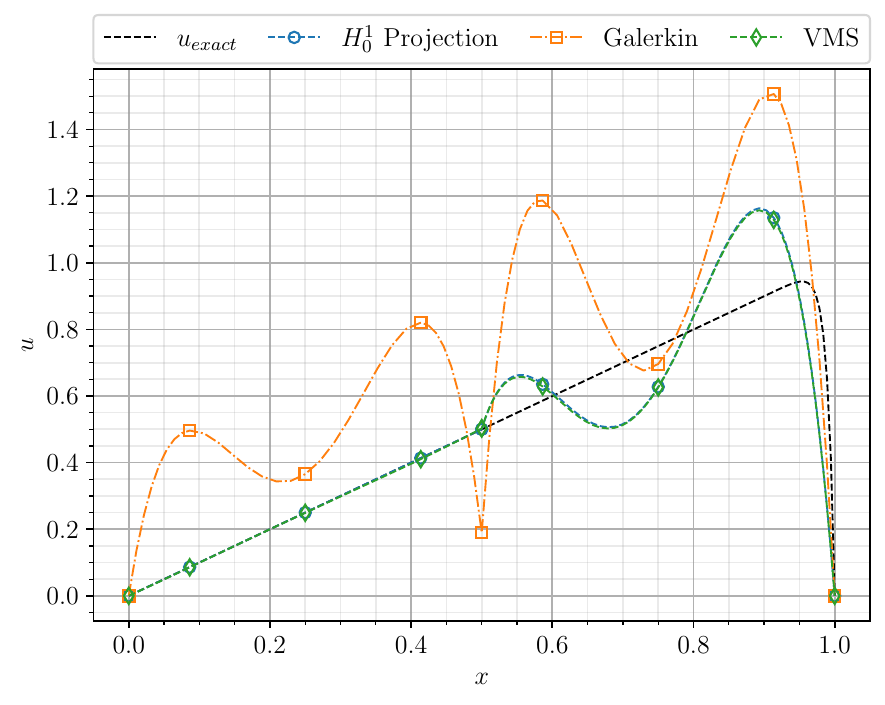}
    \caption{Polynomial degree $p = 4$ and 2 elements}
\end{subfigure}
\caption{Exact solution and its $H_0^1$ projection plotted alongside the Galerkin solution and the iterative VMS approach}
\label{fig:VMS_iter}
\end{figure}

\begin{figure}[H]
\begin{subfigure}{0.49\linewidth}
    \centering
    \includegraphics[width = \linewidth]{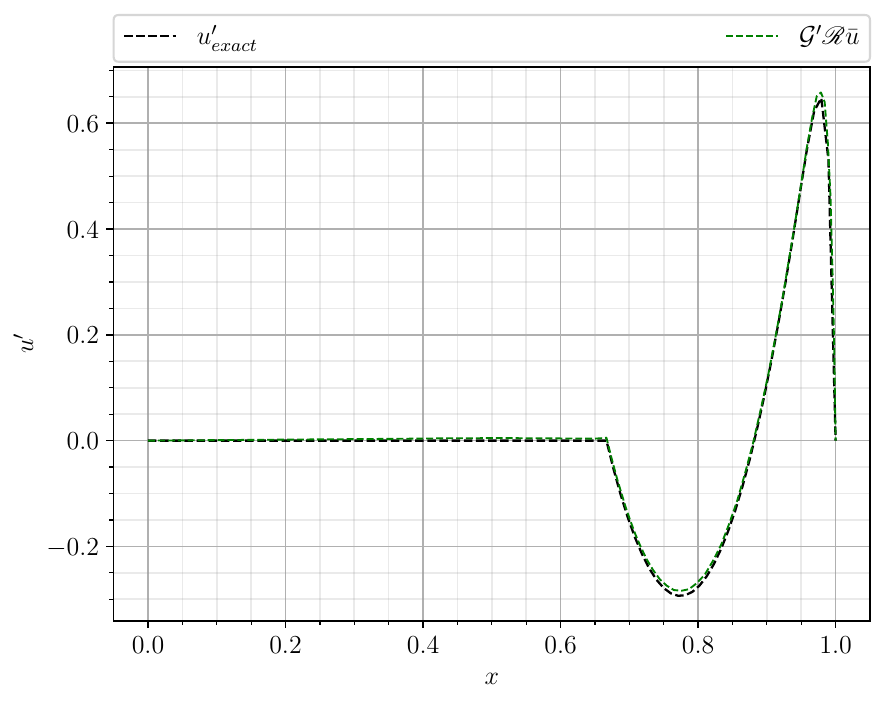}
    \caption{Polynomial degree $p = 2$ and 3 elements}
\end{subfigure}
\begin{subfigure}{0.49\linewidth}
    \centering
    \includegraphics[width = \linewidth]{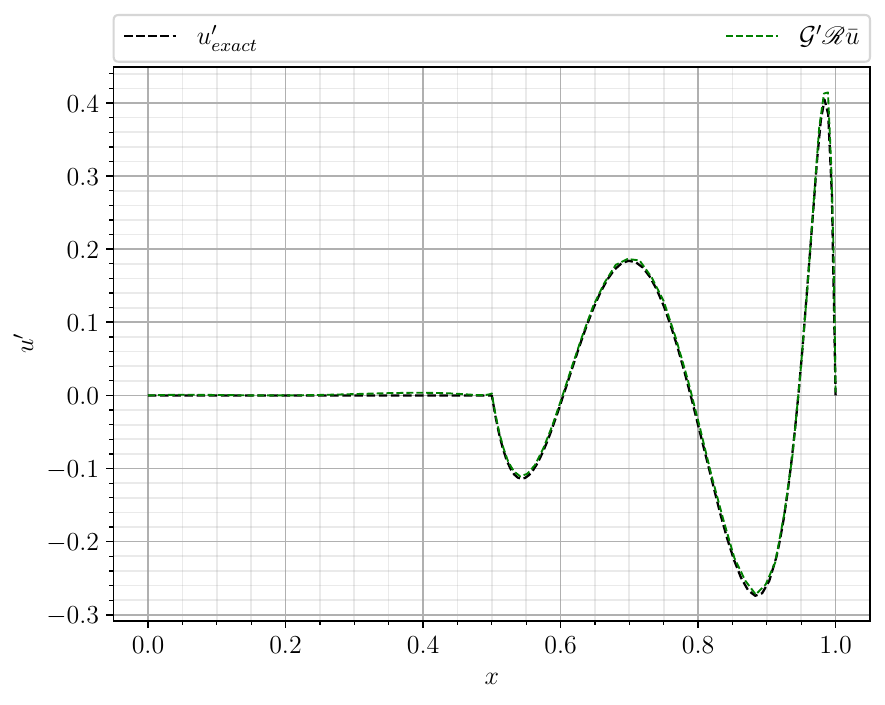}
    \caption{Polynomial degree $p = 4$ and 2 elements}
\end{subfigure}
\caption{Exact fine-scales and the fine-scales for the $H_0^1$ projection computed using the Fine-Scale Greens' function through the iterative approach}
\label{fig:VMS_iter_up}
\end{figure}

The results of this iterative VMS approach are shown in \Cref{fig:VMS_iter} and \Cref{fig:VMS_iter_up}. Firstly, \Cref{fig:VMS_iter} shows the exact solution and its $H_0^1$ projection plotted alongside the Galerkin solution and the newly proposed iterative VMS approach (labelled as 'VMS'). Subsequently, \Cref{fig:VMS_iter_up} shows the computed $u'$ alongside the exact fine-scales of the projection. From these results, it is evident that the iterative VMS approach successfully yields, up to integration and iteration error, a numerical solution that is the projection of the exact solution and it further returns the corresponding missing fine scales. It is worth noting that the proposed iterative approach is subject to some amount of numerical error as evident through the small discrepancy between the exact and the iteratively computed fine scales seen in \Cref{fig:VMS_iter_up}. This error is attributed to the choice of choice of $\varepsilon$ and the refinement level of the finer mesh where we represent $u'$ (i.e. the degree of precision of the quadrature rule used for \eqref{eq:up_ad_iter}). 

\section{Computing the Fine-Scale Greens' function for 2D problems}
\label{sec:proj_2d}
We now demonstrate that the aforementioned claim that the Fine-Scale Greens' function can be constructed using dual basis functions naturally extends to multi-dimensional cases. For the sake of demonstration, we limit our focus to the $H_0^1$ projection and the 2D Poisson equation given by
\begin{gather}
    -\nabla \cdot \nabla \phi = f, \quad \boldsymbol{x} \in \Omega = [0, 1]^2 \\
    \phi(\boldsymbol{x}) = 0, \quad \boldsymbol{x} \in \partial \Omega.
\end{gather}
The Greens' function associated with this problem may be expressed as an eigenfunction expansion as follows \cite{Haberman_2019}
\begin{gather}
    g(\boldsymbol{x}, \boldsymbol{s}) = g(x, y, s_1, s_2) =
    \displaystyle \sum_{n = 1}^{\infty} -\dfrac{(2\sin{\left(n \pi s_1 \right)} \sin{(n \pi x)})}{(n \pi \sinh{(n \pi)})} 
    \begin{cases}
        \sinh{(n \pi (s_2 - 1))} \sinh{(n \pi y)}, \quad & y < s_2 \\
        \sinh{(n \pi (y - 1))} \sinh{(n \pi s_2)}, \quad & y \geq s_2.
    \end{cases}
    \label{eq:2dPoisson_Greens}
\end{gather}
To derive the $H_0^1$ dual basis functions for the 2D case, we substitute
\begin{gather}
    \bar{\phi} := \widetilde{\boldsymbol{\psi}}^{(0)} \widetilde{\mathcal{N}}^0(\bar{\phi}) \\
    -\Delta_h v^h := \widetilde{\boldsymbol{\psi}}^{(0)} \widetilde{\mathcal{N}}^0(v^h) \quad \Longleftrightarrow \quad v^h = \left (-\Delta^{-1}_h \widetilde{\boldsymbol{\psi}}^{(0)}  \right ) \widetilde{\mathcal{N}}^0(v^h),
\end{gather}
into \eqref{eq:H01_opt}, which yields 
\begin{gather}
    \widetilde{\mathcal{N}}^0(v^h)^T \left(-\Delta^{-1}_h \widetilde{\boldsymbol{\psi}}^{(0)}, \widetilde{\boldsymbol{\psi}}^{(0)}\right)_{H^1_0(\Omega)} \widetilde{\mathcal{N}}^0(\bar{\phi}) = \widetilde{\mathcal{N}}^0(v^h)^T \left(-\Delta^{-1}_h \widetilde{\boldsymbol{\psi}}^{(0)}, \phi \right)_{H^1_0(\Omega)}, \quad \forall \widetilde{\mathcal{N}}^0(v^h)^T \in \mathbb{R}^n.
\end{gather}
As for the 1D case, we use the fact that $\left(-\Delta^{-1}_h \widetilde{\boldsymbol{\psi}}^{(0)}, \widetilde{\boldsymbol{\psi}}^{(0)}\right)_{H^1_0(\Omega)}=\delta_{i,j}$ which gives
\begin{gather}
    \widetilde{\mathcal{N}}^0(\bar{\phi}) = \left(-\Delta^{-1}_h \widetilde{\boldsymbol{\psi}}^{(0)}, \phi \right)_{H^1_0(\Omega)}.
\end{gather}
Subsequently we define our 2D $H_0^1$ dual basis functions $\boldsymbol{\mu}^{H_0^1}$ as
\begin{equation}
    \boldsymbol{\mu}^{H_0^1} := -\Delta_h^{-1} \widetilde{\boldsymbol{\psi}}^{(0)}.
\end{equation}

\begin{remark}
    The 2D $H_0^1$ dual basis functions are defined using the \textit{dual nodal} functions $\widetilde{\boldsymbol{\psi}}^{(0)}$ as opposed to the primal nodal functions used for the 1D case. See \ref{app:A} and \cite{HybridSEM} for details on the construction of the discrete Laplacian operator.
\end{remark}

\begin{figure}[H]
\begin{subfigure}{0.33\linewidth}
    \centering
    \includegraphics[width = \linewidth]{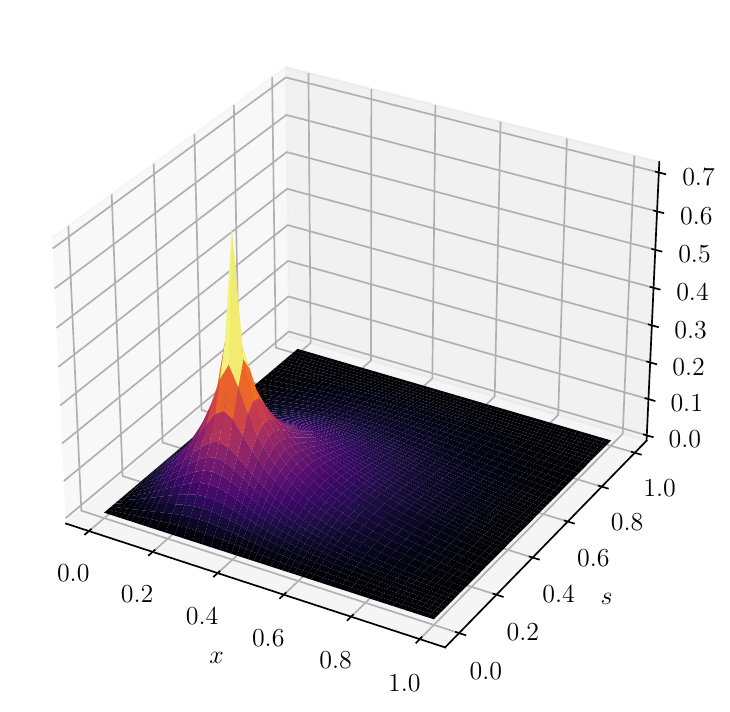}
    \caption{$(s_1, s_2) = (0.275, 0.224)$}
\end{subfigure}
\begin{subfigure}{0.33\linewidth}
    \centering
    \includegraphics[width = \linewidth]{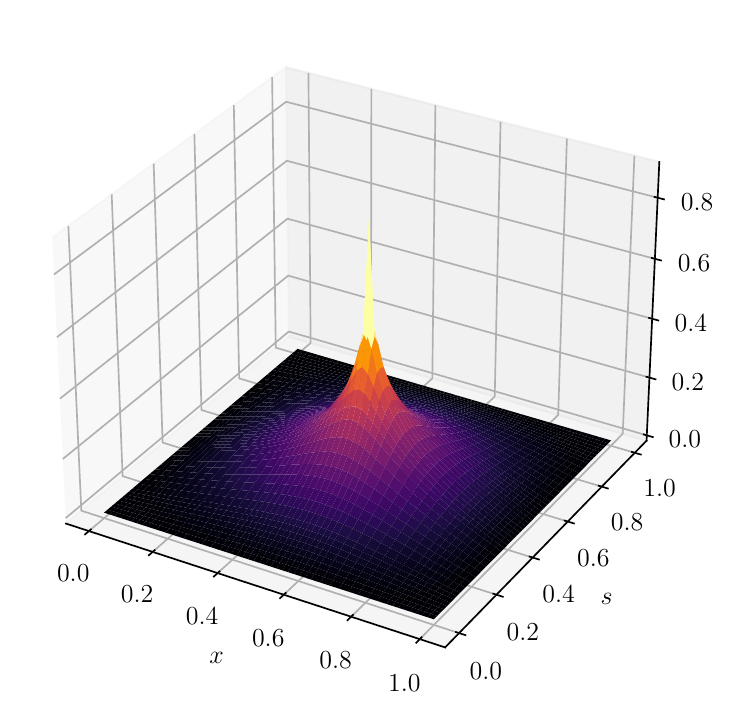}
    \caption{$(s_1, s_2) = (0.551, 0.448)$}
\end{subfigure}
\begin{subfigure}{0.33\linewidth}
    \centering
    \includegraphics[width = \linewidth]{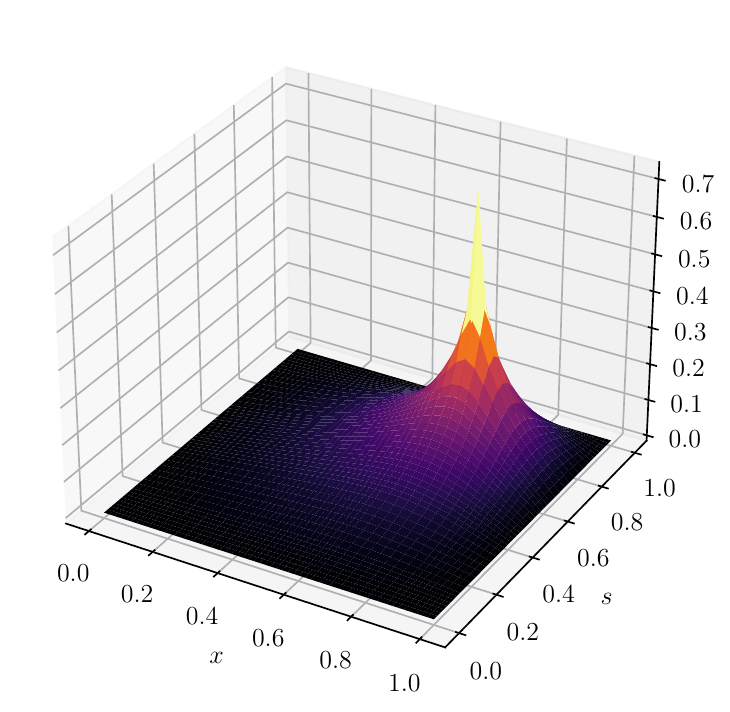}
    \caption{$(s_1, s_2) = (0.724, 0.724)$}
\end{subfigure} 
\caption{Truncated Greens' function for 2D Poisson equation at various points in the domain}
\label{fig:2d_Greens}
\end{figure}

\begin{figure}
\begin{subfigure}{0.33\linewidth}
    \centering
    \includegraphics[width = \linewidth]{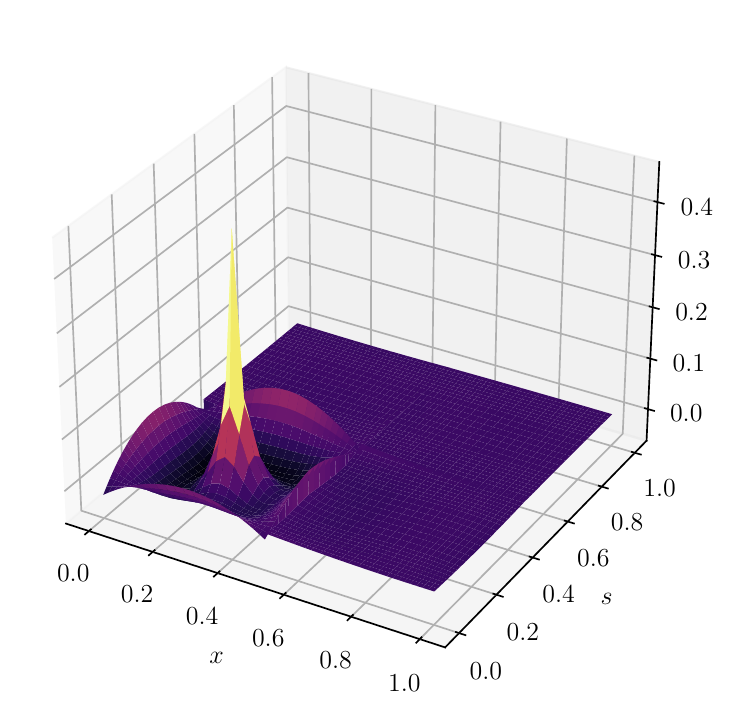}
    \caption{$(s_1, s_2) = (0.275, 0.224)$}
\end{subfigure}
\begin{subfigure}{0.33\linewidth}
    \centering
    \includegraphics[width = \linewidth]{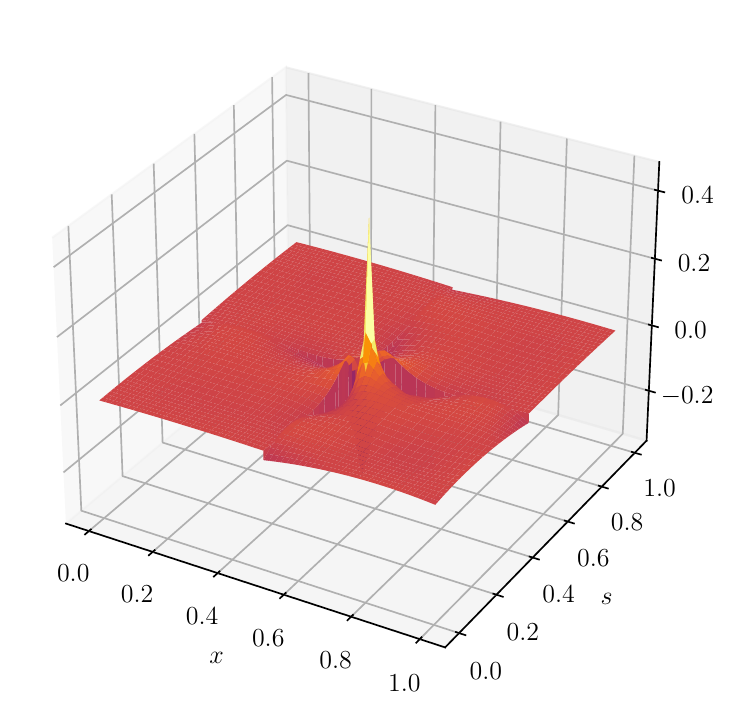}
    \caption{$(s_1, s_2) = (0.551, 0.448)$}
\end{subfigure}
\begin{subfigure}{0.33\linewidth}
    \centering
    \includegraphics[width = \linewidth]{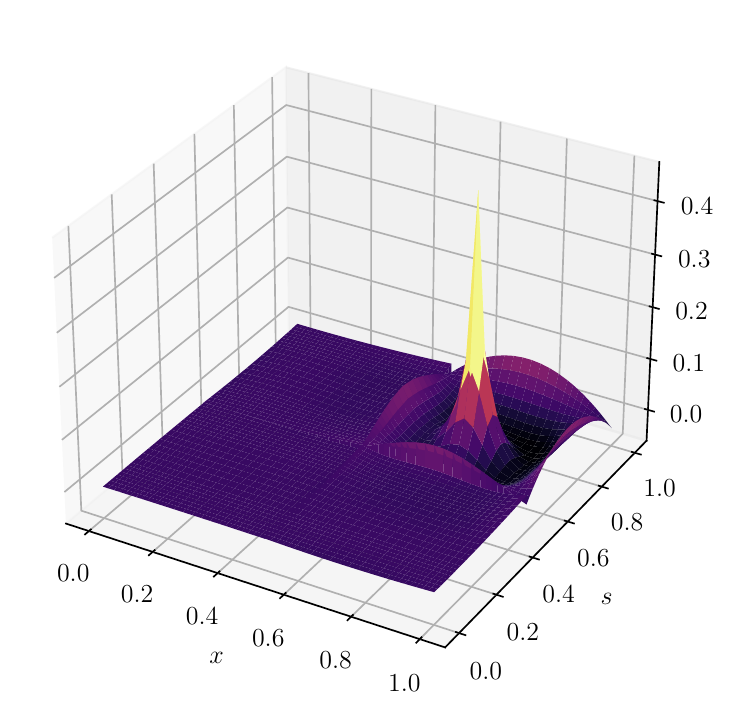}
    \caption{F$(s_1, s_2) = (0.724, 0.724)$}
\end{subfigure}
\caption{Fine-Scale Greens' function for 2D Poisson equation at various points in the domain}
\label{fig:2d_fine-Greens}
\end{figure}

Using these $H_0^1$ dual basis functions and the Greens' function in \eqref{eq:2dPoisson_Greens}, we can construct the Fine-Scale Greens' function in the exact same manner as for the 1D case. We choose to truncate the series in \eqref{eq:2dPoisson_Greens} to a finite value of 100 terms for practical implementation. The plot of the truncated Greens' function and the corresponding Fine-Scale Greens' function for the $H_0^1$ projection are visualised in \Cref{fig:2d_Greens} and \Cref{fig:2d_fine-Greens} respectively.

\subsection{Reconstruction of fine-scale terms of a projection (2D)}
\begin{figure}[H]
\begin{subfigure}{0.49\linewidth}
    \centering
    \includegraphics[width = \linewidth]{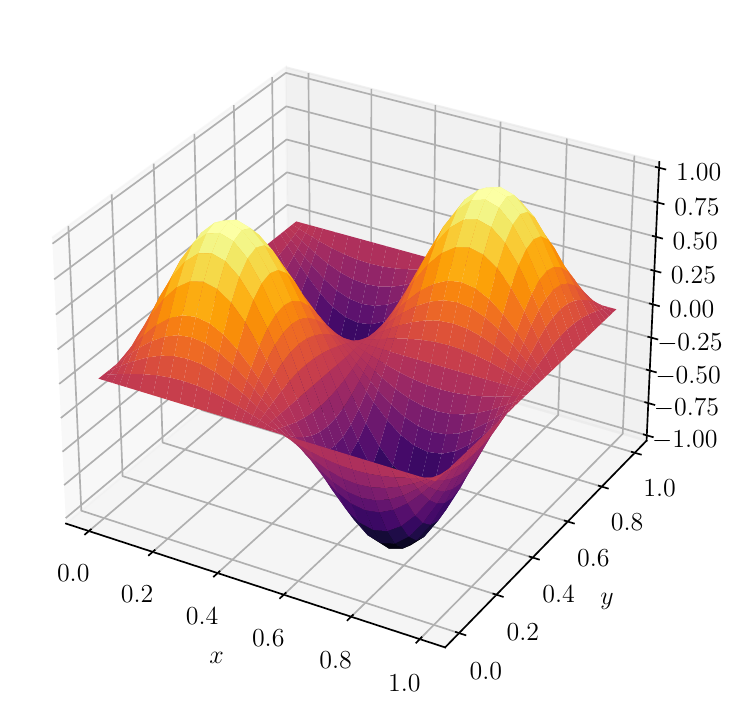}
    \caption{Exact solution}
\end{subfigure}
\begin{subfigure}{0.49\linewidth}
    \centering
    \includegraphics[width = \linewidth]{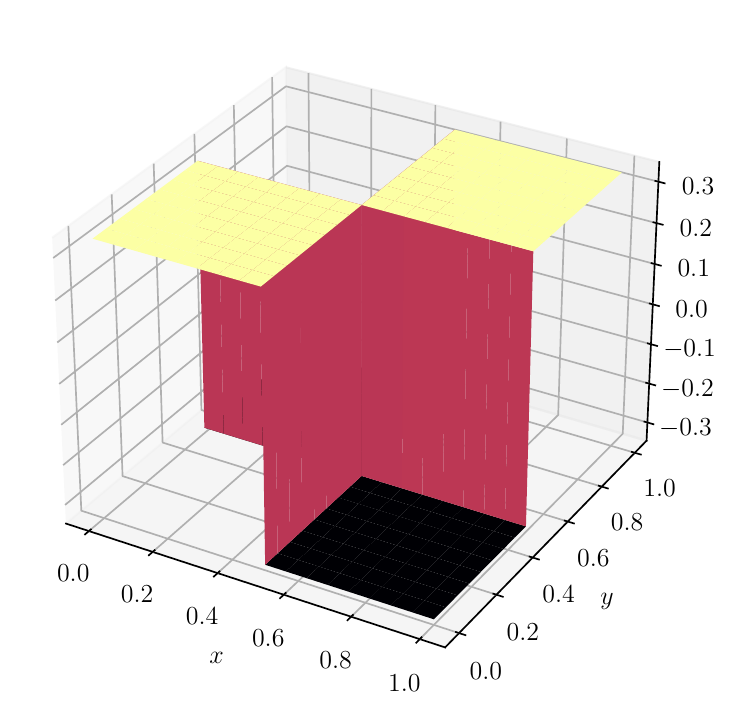}
    \caption{$H_0^1$ projection}
\end{subfigure}
\caption{Exact solution and $H_0^1$ projection of 2D Poisson equation on a mesh with $2 \times 2$ elements with polynomial degree $p = 1$}
\label{fig:u_ex_2d_proj_p1}
\end{figure}

\begin{figure}[H]
\begin{subfigure}{0.49\linewidth}
    \centering
    \includegraphics[width = \linewidth]{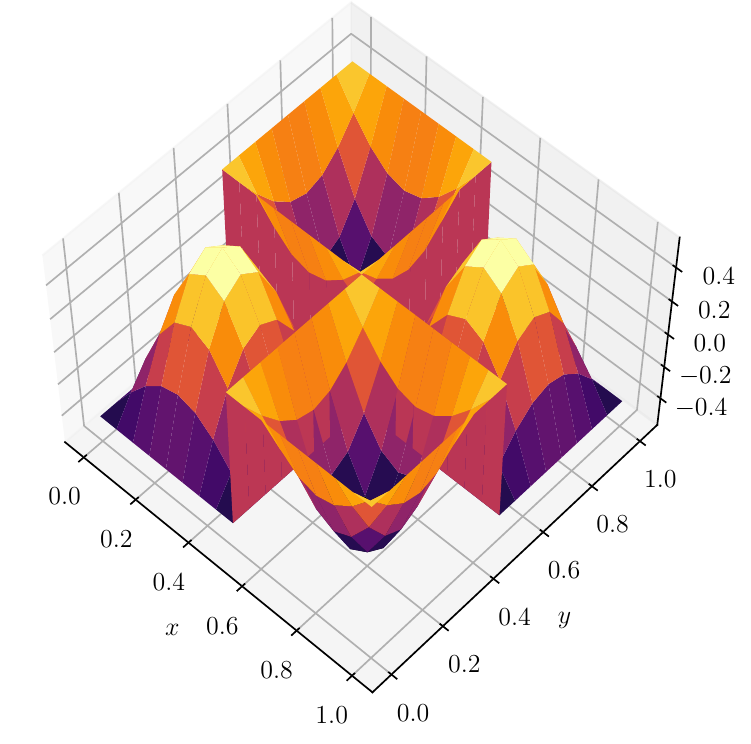}
    \caption{Exact fine-scales for the $H_0^1$ projection}
\end{subfigure}
\begin{subfigure}{0.49\linewidth}
    \centering
    \includegraphics[width = \linewidth]{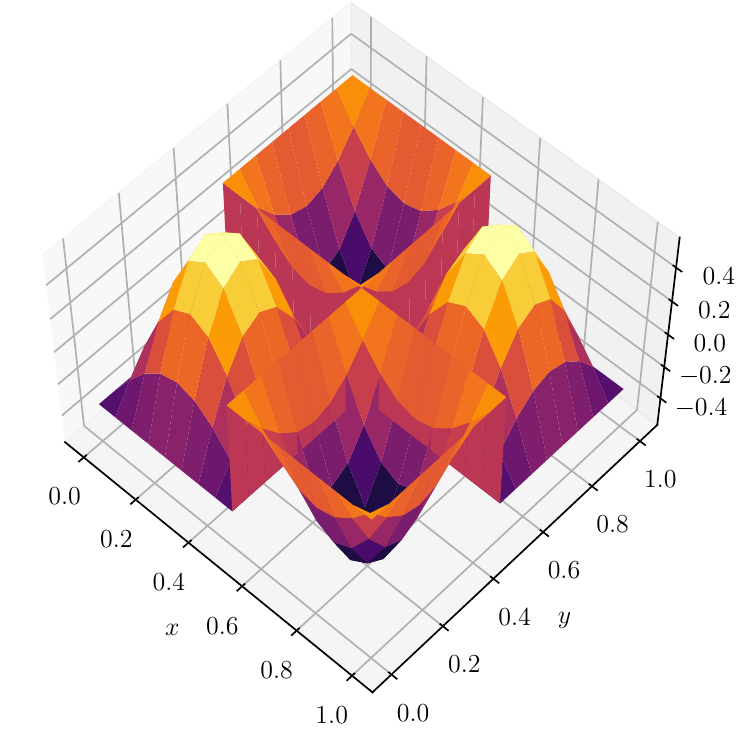}
    \caption{Fine-scales computed using the Fine-Scale Greens' function}
\end{subfigure}
\caption{Unresolved fine-scales of the $H_0^1$ projection of the 2D Poisson equation on a mesh with $2 \times 2$ elements polynomial degree $p = 1$}
\label{fig:u_p_2d_proj_p1}
\end{figure}

\begin{figure}[H]
\begin{subfigure}{0.49\linewidth}
    \centering
    \includegraphics[width = \linewidth]{Images/u_exact.pdf}
    \caption{Exact solution}
\end{subfigure}
\begin{subfigure}{0.49\linewidth}
    \centering
    \includegraphics[width = \linewidth]{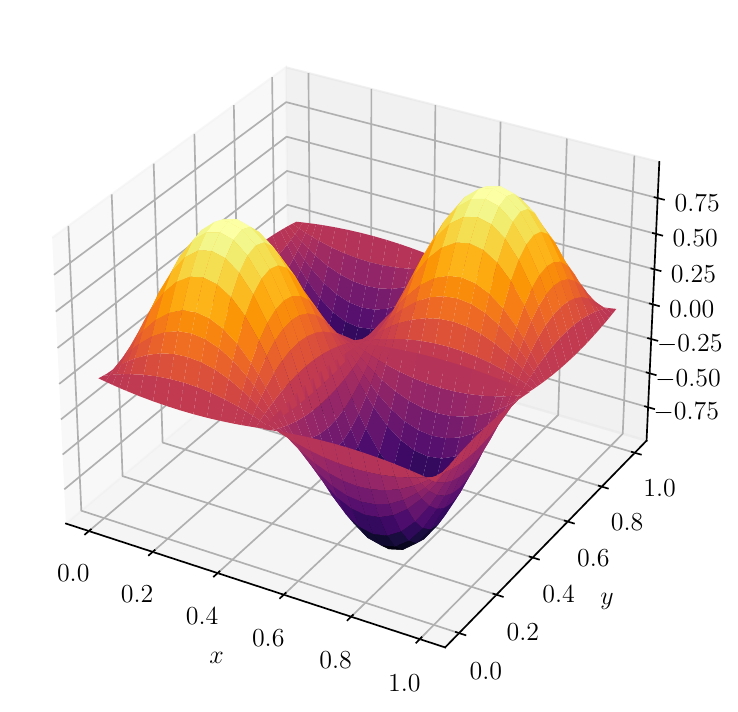}
    \caption{$H_0^1$ projection}
\end{subfigure}
\caption{Exact solution and $H_0^1$ projection of 2D Poisson equation on a mesh with $2 \times 2$ elements with polynomial degree $p = 4$}
\label{fig:u_ex_2d_proj_p2}
\end{figure}

\begin{figure}[H]
\begin{subfigure}{0.49\linewidth}
    \centering
    \includegraphics[width = \linewidth]{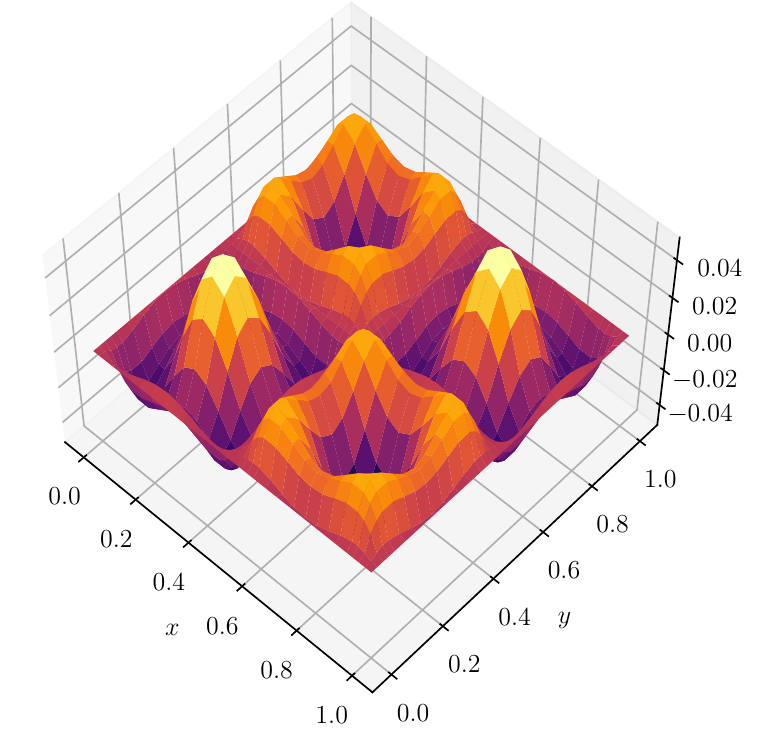}
    \caption{Exact fine-scales for the $H_0^1$ projection}
\end{subfigure}
\begin{subfigure}{0.49\linewidth}
    \centering
    \includegraphics[width = \linewidth]{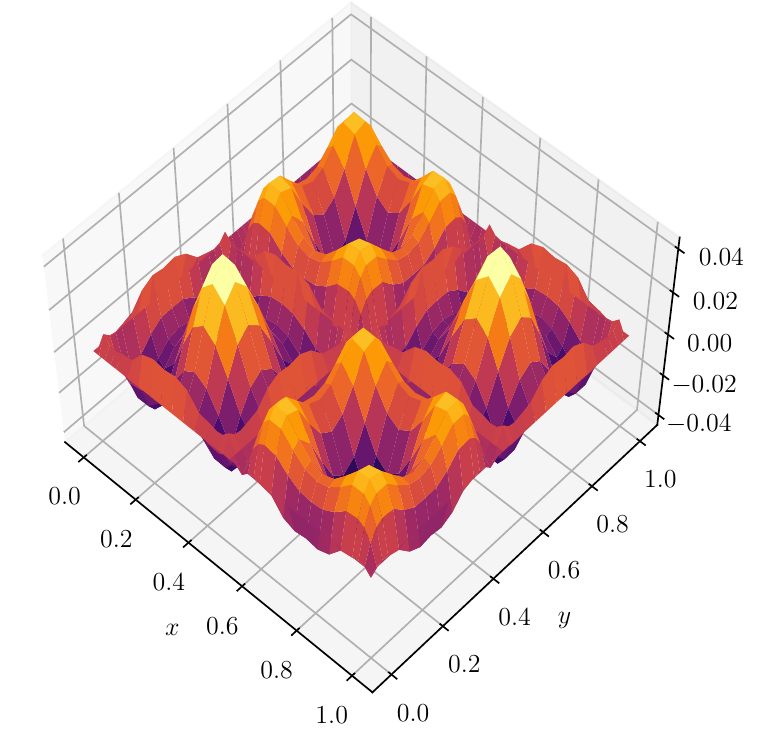}
    \caption{Fine-scales computed using the Fine-Scale Greens' function}
\end{subfigure}
\caption{Unresolved fine-scales of the $H_0^1$ projection of the 2D Poisson equation on a mesh with $2 \times 2$ elements with polynomial degree $p = 4$}
\label{fig:u_p_2d_proj_p2}
\end{figure}

We perform similar numerical tests as for the 1D case where we compute the projection of an exact solution to the Poisson problem and reconstruct the fine scales. For this 2D case, we take as the source term
\begin{equation}
    f = 8 \pi^2 \sin{(2 \pi x)} \sin{(2 \pi y)},
\end{equation}
which leads to the following exact solution
\begin{equation}
    \phi_{exact} = \sin{(2 \pi x)} \sin{(2 \pi y)}.
\end{equation}
The exact solution and its $H_0^1$ projection onto two different meshes\footnote{$p$ refers to degree of nodal polynomials, see \Cref{rem:edge_p}} are shown in \Cref{fig:u_ex_2d_proj_p1} and \Cref{fig:u_ex_2d_proj_p2}. Subsequently, the exact fine scales and the reconstruction thereof computed using the Fine-Scale Greens' function are shown in \Cref{fig:u_p_2d_proj_p1} and \Cref{fig:u_p_2d_proj_p2}. These plots clearly indicate that the described Fine-Scale Greens' function correctly reproduces the missing unresolved scales of a projection. 

%% file: sections/06-summary.tex
\section{Summary}
\label{sec:summary}


In this paper, we have proposed a new approach for explicitly constructing the Fine-Scale Greens' function by employing the concept of dual basis functions. We have shown that these dual basis functions encode the projector and are explicitly computable thus making them an appropriate choice for the $\boldsymbol{\mu}$'s presented in \cite{hughes_sangalli_2007}. Moreover, the constructed Fine-Scale Greens' function is clearly shown to successfully reconstruct all the missing fine scales truncated with the projection. Since the dual basis functions naturally extend to higher dimensions, we have a generalised approach for computing the Fine-Scale Greens' function for any projector in arbitrary dimensions. Furthermore, even though we solely focus on the Fine-Scale Greens' functions for the Poisson equation, we have demonstrated that it may be employed in VMS approaches for other problems such as the steady advection-diffusion problem albeit in a non-monolithic manner. We do remark that while everything we presented was implemented in the context of the MSEM, the concepts of the dual basis may be employed in any other framework to produce an appropriate Fine-Scale Greens' function. Regarding future work, a natural continuation of this paper will be the development of numerical integration rules that accurately and efficiently integrate the Fine-Scale Greens' function. Subsequently, the methodology could be expanded to address increasingly intricate non-linear problems, such as the Navier-Stokes equations, with consideration given to multi-dimensional advection-diffusion problems as intermediary stages.

%% file: sections/08-appendix.tex
\section{Discrete Laplacian operator in 2D}
\label{app:A}
In this section, we shall describe the construction of the discrete Laplacian operator used for defining the $H_0^1$ projection in 2D from \Cref{sec:proj_2d}. We start with a Poisson problem defined as follows
\begin{gather}
    -\nabla \cdot \nabla \phi = f, \quad \boldsymbol{x} \in \Omega \\
    \phi(\boldsymbol{x}) = 0, \quad \boldsymbol{x} \in \partial \Omega.
\end{gather}
We then rewrite the problem in a mixed formulation by defining a vector field $\underline{u}$ as the gradient of $\phi$
\begin{gather}
    \underline{u} - \nabla \phi = 0 \label{eq:mixed1}\\
    -\nabla \cdot \underline{u} = f. \label{eq:mixed2}
\end{gather}
The Galerkin weak form of this mixed formulation is obtained by testing \eqref{eq:mixed1} with $\underline{v} \in H(div, \Omega)$ and \eqref{eq:mixed2} with $\varphi \in L^2(\Omega)$
\begin{gather}
    \left(\underline{v}, \underline{u} \right)_{L^2(\Omega)} + \left(\nabla \cdot \underline{v}, \phi \right)_{L^2(\Omega)} = \int_{\partial \Omega} \phi \underline{v} \cdot \underline{n} \: \mathrm{d} \Gamma, \quad \forall \underline{v} \in H(div, \Omega) \label{eq:wk_1} \\
    -\left(\varphi, \nabla \cdot \underline{u} \right)_{L^2(\Omega)} = \left(\varphi, f \right)_{L^2(\Omega)}, \quad \forall \varphi \in L^2(\Omega), \label{eq:wk_2}
\end{gather}
where we have performed integration by parts on the gradient of $\phi$ and the emerging boundary integral is used to weakly impose the homogeneous Dirichlet boundary condition by setting it to zero. Carrying out this integration by parts is what allows us to eventually perform the $H_0^1$ projection in a discontinuous space.

An important remark about the choice of function spaces is that they satisfy the De Rham sequence, i.e. applying the divergence to an element of $H(div, \Omega)$ maps the element to $L^2(\Omega)$
\begin{equation*}
    H(div, \Omega) \xrightarrow{\nabla \cdot} L^2(\Omega).
\end{equation*} 
This property must also be satisfied at the discrete level when selecting the finite-element subspaces, see \cite{Arnold_Falk_Winther_2010}. 


Transforming \eqref{eq:wk_1} and \eqref{eq:wk_2} into the discrete setting by expanding $\phi$ in the space of dual nodal basis functions ($\phi = \widetilde{\boldsymbol{\psi}}^{(0)} \boldsymbol{\phi}$) and $\underline{u}$ in (2D) edge basis functions ($\underline{u} = \boldsymbol{\psi}^{(1)} \boldsymbol{u}$) gives
\begin{gather}
    \MassM{1} \boldsymbol{u} + \nabla^T_h \cdot \boldsymbol{\phi} = 0 \label{eq:disc_wk1} \\
    -\nabla_h \cdot \boldsymbol{u} = \boldsymbol{\hat{f}}, \label{eq:disc_wk2} 
\end{gather}
where $\MassM{1}$ is mass matrix whose entries correspond to the $L^2$ inner product of the basis functions spanning a subset of $H(div, \Omega)$ (2D edge basis functions $\left(\psi^{(1)}_i, \psi^{(1)}_j \right)_{L^2(\Omega)}$), $\nabla_h \cdot$ is the \textit{discrete} divergence operator, and $\boldsymbol{u}$ and $\boldsymbol{\phi}$ are the degrees of freedom. Given the invertibility of $\MassM{1}$, \cite{jain_zhang_palha_gerritsma_2021}, we rearrange \eqref{eq:disc_wk1} to get an expression for $\boldsymbol{u}$ as follows
\begin{equation}
    \boldsymbol{u} = -\MassMinv{1} \nabla^T_h \cdot \boldsymbol{\phi},
\end{equation}
which we plug into \eqref{eq:disc_wk2} to get
\begin{equation}
    \nabla_h \cdot \MassMinv{1} \nabla^T_h \cdot \boldsymbol{\phi} = \boldsymbol{\hat{f}}.
\end{equation}
Hence, the discrete operator for the negative Laplacian is given as
\begin{equation}
    -\Delta_h = \nabla_h \cdot \MassMinv{1} \nabla^T_h \cdot.
\end{equation}